\documentclass[a4paper,12pt,oneside,reqno]{amsart}

\usepackage[T1]{fontenc}
\usepackage[utf8]{inputenc}

\usepackage[a4paper]{geometry}
\usepackage{amssymb,amsmath,amsthm}
\usepackage{times}

\usepackage{mathtools}
\usepackage{cite}
\usepackage[mathcal]{euscript}
\usepackage{hyperref}
\hypersetup{hidelinks,hypertexnames=false}
\newenvironment{Proof}[1][Proof.]{\begin{proof}[\upshape\bfseries #1]}{\end{proof}}
\usepackage{xcolor}

\usepackage{graphicx}

\def\dim{\mathrm{dim}}

\def\As{\mathrm{As}}

\def\id{\mathop{\fam 0 id}\nolimits}

\let\<\langle
\let\>\rangle

\newcommand{\kk}{\mathbb{k}}
\newcommand{\NB}{\mathrm{noncom\text{-}Bicom}}

\newcommand{\AsS}{\mathrm{As}(S)}
\newcommand{\Ker}{\operatorname{Ker}}

\newcommand{\Ical}{\mathcal I}
\newcommand{\Ncal}{\mathcal N}
\newcommand{\Xcal}{\mathcal X}
\newcommand{\Ycal}{\mathcal Y}

\newcommand{\LL}{\mathsf L}
\newcommand{\RR}{\mathsf R}
\newcommand{\arity}{\operatorname{ar}}
\newcommand{\Mag}{\operatorname{Mag}_{\succ}}
\newcommand{\Id}{\operatorname{Id}}
\newcommand{\Ccal}{\mathcal C}

\newcommand{\outdeg}{\operatorname{outdeg}}
\newcommand{\Span}{\operatorname{span}_{\kk}}

\newtheorem{definition}{Definition}
\newtheorem{lemma}{Lemma}
\newtheorem{proposition}{Proposition}
\newtheorem{theorem}{Theorem}
\newtheorem{corollary}{Corollary}
\newtheorem{remark}{Remark}

\allowdisplaybreaks[1]

\title[Special identities of nonsymmetric bicommutative operations]{Special identities of nonsymmetric bicommutative operations in associative algebras with an operator}
\author{Bauyrzhan Sartayev$^{*}$}
\address{Narxoz University, Almaty, Kazakhstan}
\address{Institute of Mathematics and Mathematical Modeling, Almaty, Kazakhstan}
\email{baurjai@gmail.com}
\author{Abdibek Ydyrys}
\address{SDU University, Kaskelen, Kazakhstan}
\email{abdibek54@gmail.com}
\keywords{Nonsymmetric operad, bicommutative algebra, special identity, associative algebra with an operator}
\subjclass[2020]{17A30, 17A50, 16R10}
\thanks{${}^{*}$Corresponding author: Bauyrzhan Sartayev (baurjai@gmail.com)}

\begin{document}
\begin{abstract}
We study the identities of the operations $a\prec b=aS(b)$ and $a\succ b=S(a)b$ in associative algebras with a linear operator satisfying $aS(S(b)c)=S(aS(b))c$. Denote the corresponding operated operad by $\As^{(S)}$. We obtain a complete minimal generating set of special identities and prove that this set is infinite. More precisely, in every arity $n\geq8$, there are exactly $n-2$ independent new identities which do not follow from identities of smaller arity. Thus the natural morphism 
\[
\phi\colon\NB\to\As^{(S)}
\]
has an infinitely generated kernel, although its source is a finitely presented quadratic operad. To the best of our knowledge, this is the first example of this phenomenon for a derived realization of a quadratic operad.
\end{abstract}
\maketitle

\section{Introduction}
An algebra with a linear operator often gives rise to new algebraic
structures on the same vector space. The identities of the new
operations depend both on the original multiplication and on the
identity satisfied by the operator. A natural question is to determine
all these identities and to decide whether a finite list is sufficient.
In the language of operads, this amounts to describing the kernel of
the morphism defined by the new operations; see
\cite{GK,LV} for the general operadic setting.

For a derivation $d$, the operations $a\prec b=ad(b)$ and
$a\succ b=d(a)b$ have been studied through derived varieties and
Manin products in \cite{KSO,Vallette}. The commutative case
is related to Novikov algebras and their free algebras
\cite{DL}. The associative case leads to noncommutative Novikov
algebras, for which embeddings into differential associative algebras
and a monomial basis were obtained in \cite{SKNov,DSBasis}.
Related embedding questions for differential Perm-algebras and
Zinbiel algebras were considered in \cite{KSPerm,KMS}. Special
identities also occur for the operations obtained from differential
Poisson algebras; the Gelfand--Dorfman case was studied in
\cite{KSO,KSGD}.

Another source of such constructions is a Rota--Baxter operator.
The products $aR(b)$ and $R(a)b$, for an operator of weight zero,
satisfy the dendriform identities of Loday \cite{Loday}.
Their relation to free Rota--Baxter algebras and enveloping algebras
was studied in \cite{EFGuo}; see also \cite{Guo}.
Aguiar \cite{Aguiar} considered analogous constructions for
pre-Poisson algebras. These examples suggest studying the identities
of $aS(b)$ and $S(a)b$ for other relations on a linear operator $S$.

In this paper, we consider an associative non-unital algebra $A$
with a linear operator $S$ satisfying
\[
aS(S(b)c)=S(aS(b))c.
\]
Then $a\prec b=aS(b)$ and $a\succ b=S(a)b$ satisfy
\[
a\succ(b\prec c)=(a\succ b)\prec c,
\qquad a\prec(b\succ c)=(a\prec b)\succ c.
\]
These identities define the nonsymmetric bicommutative operad
$\NB$ considered in \cite{MS}. The same paper raises the question
whether the corresponding morphism into the operad of associative
algebras with $S$ is injective. We give a negative answer and
describe its entire kernel.

The first special identity already occurs in arity four:
\[
(a\prec b)\succ(c\succ d)=((a\prec b)\prec c)\succ d.
\]
Its two sides are distinct normal monomials in $\NB$, although their
operated images are equal. Thus the cubic relations do not describe
all identities of the derived operations. One must also determine
whether the identities obtained in higher arities follow from the
earlier ones. It is this second question that leads to the quotient
$Q_n$ considered below.

Recall that a bicommutative algebra satisfies
$(ab)c=(ac)b$ and $a(bc)=b(ac)$. Free bicommutative algebras were
described in \cite{DIT}. Their polynomial identities and subvarieties
were studied in \cite{DI,DrenskyBicom}, and the finite basis problem
for their subvarieties was solved in \cite{DZ} over an arbitrary
field. Here we retain two operations and work with nonsymmetric
substitutions. In this setting the ideal of special identities is
not finitely generated.

Let $\phi\colon\NB\to\AsS$ be the morphism given by the derived
operations. Put $K_n=\ker\phi_n$, and let $D_n$ be the space of
arity-$n$ consequences of the components $K_j$ with $j<n$.
Thus $Q_n=K_n/D_n$ records the identities which cannot be obtained
from smaller arities. We introduce explicit families $F_r,G_r,J_r$
and $H_{p,q}$ and prove that, for every $n\geq8$, the classes of
$$F_{n-2},\ \rho F_{n-2},\ G_{n-3},\ \rho G_{n-3},\
J_{n-4},\ \rho J_{n-4},H_{p,n-p-1},\qquad\qquad 4\leq p\leq n-5,$$
form a basis of $Q_n$, where $\rho$ reverses the variables and
interchanges $\prec$ and $\succ$. Consequently,
$\dim Q_n=n-2$. We also give the initial generators in arities
$4,5,6,7$, whose numbers are $1,4,7,6$, respectively. Together
these identities form a minimal homogeneous generating set of
$\ker\phi$. All the statements hold over an arbitrary field.

We also determine all identities involving only one of the two
operations. Since the defining relations of $\NB$ are mixed, the
suboperad generated by $\succ$ alone is free. Its map into $\AsS$
nevertheless has a nonzero kernel. We prove that this kernel is
generated by $F_r,G_r$, $r\geq3$, and two additional identities
$K_6,E_6$ of arity six. Reflection gives the corresponding result
for $\prec$. This description is used in the proof for the two
operations together.

The proof uses normal monomials and operated words. Rewriting methods
for associative algebras and operads are developed in
\cite{Bergman,DK,BD}; their versions for algebras with operators
are considered in \cite{BCQ,LM}. For the relation studied here,
each occurrence of $S$ determines an interval of variables, and
the set of singleton intervals is preserved by every transformation.
This invariant separates the kernel into smaller subspaces.
We describe normal forms for monomials containing only $\succ$,
study the subspace with singleton set $\{1,2,n-1,n\}$, and reduce
the remaining cases to these two. Coefficient functionals then
prove that the resulting generators are independent modulo
consequences of smaller arity.

A point which requires care is that a transformation between two
derived monomials can pass through operated words which are not
themselves derived monomials. The literal images are the regular
words: at every concatenation level there is exactly one variable
outside the arguments of $S$. We describe their binary preimages
and prove the contraction statements needed to pass from an
operated-word equivalence to a consequence of a smaller binary
identity. These statements allow us to use transformations in
either direction throughout the proof.

Finally, we consider the case when the associative multiplication
is commutative. An explicit polynomial construction embeds the free
bicommutative algebra into an algebra with $S$ under the product
$a\circ b=aS(b)$. Hence there are no additional special identities
in this case. The construction already works with $S^2=S$ and
$S|_{A^2}=\id$.

\section{Preliminaries and special identities}
We work over a field $\kk$. All associative algebras are non-unital, and
all substitution blocks are nonempty. The variables of a multilinear
monomial occur in the order $x_1,\ldots,x_n$, and all operadic
substitutions are nonsymmetric.

Let $\AsS$ be the operad of associative algebras with a linear operator
$S$ satisfying
\begin{equation}
A S(S(B)C)=S(A S(B))C.
\label{R}
\end{equation}
We write $u\sim_S v$ for congruence modulo this relation. Define
\begin{equation}
a\prec b=aS(b),\qquad a\succ b=S(a)b.
\label{derived}
\end{equation}
Since both $a\succ(b\prec c)$ and $(a\succ b)\prec c$ equal $S(a)bS(c)$,
and the second identity below is \eqref{R}, we obtain
\begin{equation}
a\succ(b\prec c)=(a\succ b)\prec c,
\qquad a\prec(b\succ c)=(a\prec b)\succ c.
\label{NB}
\end{equation}
These identities define $\NB$, and \eqref{derived} defines
$\phi\colon\NB\to\AsS$. We retain the notation $K_n=\ker\phi_n$,
\begin{equation*}
D_n=\langle K_1\oplus\cdots\oplus K_{n-1}\rangle_{\mathrm{op}}(n),
\qquad Q_n=K_n/D_n.
\end{equation*}

\noindent\textbf{Normal monomials and the kernel.}

Write $x(u,v)=u\prec v$ and $y(u,v)=u\succ v$. Let $\Xcal$ and $\Ycal$
be the classes of monomials containing only $\prec$ and only $\succ$,
respectively; we call these monomials pure. A single leaf, denoted here
by $1$, belongs to both classes.
Thus
\begin{equation*}
\Xcal=\{1\}\sqcup x(\Xcal,\Xcal),\qquad
\Ycal=\{1\}\sqcup y(\Ycal,\Ycal).
\end{equation*}
Define
\begin{equation*}
\Ncal=\{1\}\sqcup x(\Ncal,\Xcal)\sqcup y(\Ncal,\Ycal).
\end{equation*}
We call these monomials normal. Every nontrivial normal monomial has the
form $U\prec P$ or $U\succ Q$, where $U\in\Ncal$, $P\in\Xcal$ and
$Q\in\Ycal$.

\begin{lemma}
The normal monomials form a basis of $\NB$.
\end{lemma}
\begin{Proof}
Let $\alpha,\beta$ be the two distinct operations. The defining relations
have the form
\[\alpha(A,\beta(B,C))=\beta(\alpha(A,B),C).\]
On the vector space with basis $\Ncal$, define $A*_{\alpha}B$ recursively.
If $B$ contains only $\alpha$, put $A*_{\alpha}B=\alpha(A,B)$.
Otherwise, write $B=\gamma(B_L,B_R)$ in normal form and put
\[
A*_{\alpha}B=\gamma(A*_{\alpha}B_L,B_R).
\]
The right operand has smaller arity at each step. Hence the product is
well defined, is normal, and preserves the number of vertices of each type.

We verify the defining identities by induction on the arity of $C$.
If $C$ contains only $\beta$, the definition gives
$A*_{\alpha}(B*_{\beta}C)=(A*_{\alpha}B)*_{\beta}C$.
Otherwise, write $C=\gamma(C_L,C_R)$ and put $D=B*_{\beta}C_L$.
Since $D$ contains $\beta$, the monomial $\gamma(D,C_R)$ is not pure
$\alpha$. Therefore
\[
A*_{\alpha}(B*_{\beta}C)
=\gamma(A*_{\alpha}D,C_R)
=\gamma((A*_{\alpha}B)*_{\beta}C_L,C_R)
=(A*_{\alpha}B)*_{\beta}C.
\]

The same recursion reduces monomials modulo \eqref{NB}. When
$\gamma\ne\alpha$, its root step is a defining relation. When
$\gamma=\alpha$, the factor $B_L$ contains $\beta$. Applying the mixed
relations along the path from a nearest $\beta$-vertex to the root gives
$B_L=\beta(P,Q)$. Then
\begin{multline*}
\alpha(A,\alpha(\beta(P,Q),C))
=\alpha(A,\beta(P,\alpha(Q,C)))
=\beta(\alpha(A,P),\alpha(Q,C))\\
=\alpha(\beta(\alpha(A,P),Q),C)
=\alpha(\alpha(A,\beta(P,Q)),C).
\end{multline*}
Thus every monomial is equal to a normal one. Evaluating a normal monomial
in the algebra just constructed gives the corresponding basis vector.
Consequently, normal monomials are linearly independent in $\NB$.
\end{Proof}

Let $\Ncal_n$ be the normal basis in arity $n$. Relation \eqref{R}
identifies operated words in pairs. Hence the quotient of the free
operated associative algebra has a basis consisting of their congruence
classes. For $T\in\Ncal_n$, put $\nu_n(T)=[\phi(T)]_S$.

\begin{lemma}
\label{lem:fibres}
Let $\Ncal_n=C_1\sqcup\cdots\sqcup C_r$ be the fibres of $\nu_n$, and
choose $T_i\in C_i$. Then
\[
\{T-T_i:T\in C_i\setminus\{T_i\},\ 1\leq i\leq r\}
\]
is a basis of $K_n$.
\end{lemma}
\begin{Proof}
Since distinct operated-word classes are linearly independent,
$\phi(\sum_T\alpha_TT)=0$ if and only if
$\sum_{T\in C_i}\alpha_T=0$ for every $i$. The subspace with this
condition on $C_i$ has basis $T-T_i$, $T\in C_i\setminus\{T_i\}$.
Taking the direct sum over the fibres proves the assertion.
\end{Proof}

\noindent\textbf{The special identities.}

For $r\geq1$, let $L_r$ be the left $\succ$-comb, with $L_1(u)=u$:
\begin{equation*}
L_r(u_1,\ldots,u_r)=(((u_1\succ u_2)\succ u_3)\cdots)\succ u_r.
\end{equation*}
Put
\begin{align*}
T_4(x_1,x_2,x_3,x_4)&=x_1\succ((x_2\succ x_3)\succ x_4),\\
T_{k+1}(x_1,\ldots,x_{k+1})&=T_k(x_1,\ldots,x_k)\succ x_{k+1},
\qquad k\geq4.
\end{align*}
Reflection is defined by
\begin{equation*}
\rho(u\succ v)=\rho(v)\prec\rho(u),\qquad
\rho(u\prec v)=\rho(v)\succ\rho(u),
\end{equation*}
followed by renumbering the leaves from left to right. Write
$Q_k=\widetilde T_k=\rho(T_k)$. Reflection preserves \eqref{R}.

We use the following identities, together with their reflections:
\begin{align*}
F_r:\quad a\succ\bigl(b\succ L_r(c_1,\ldots,c_r)\bigr)
&=a\succ L_{r-1}\bigl(b\succ(c_1\succ c_2),c_3,\ldots,c_r\bigr),
\quad r\geq3,\\
G_r:\quad a\succ\bigl(L_r(b_1,\ldots,b_r)\succ(c\succ d)\bigr)
&=T_{r+1}(a,b_1,\ldots,b_r)\succ(c\succ d),
\quad r\geq3,\\
J_r:\quad (a\succ(b\succ c))\prec(d\prec\widetilde T_r)
&=\bigl(a\succ((b\succ c)\succ d)\bigr)\prec\widetilde T_r,
\quad r\geq4,\\
H_{p,q}:\quad T_p\prec(z\prec\widetilde T_q)
&=(T_p\succ z)\prec\widetilde T_q,
\quad p,q\geq4.
\end{align*}
Their arities are $r+2$, $r+3$, $r+4$ and $p+q+1$, respectively.

\noindent\textbf{Validity of the identities.}
\label{sec:validity}

\begin{definition}
For an operated word $u$, write $\LL(u)$ if its $S$-class contains a word
$AS(B)$ with $A,B$ nonempty. Similarly, write $\RR(u)$ if its class
contains a word $S(A)B$ with $A,B$ nonempty.
\end{definition}

\begin{lemma}
\label{lem:transport}
If $\LL(u)$, then
\begin{align}
xS(u)d&=S(xu)d,\label{K0}\\
S(u)S(v)d&=S(uS(v))d,\label{H0}\\
S(ucS(d))e&=S(u)cS(d)e.\label{J0}
\end{align}
\end{lemma}
\begin{Proof}
Choose $u=AS(B)$ in its congruence class. For \eqref{K0}, both sides
equal $xA S(S(B)d)$. For the other two identities, apply \eqref{R}:
\[
S(AS(B)S(v))d=AS(B)S(S(v)d)
=AS(S(B)S(v))d=S(AS(B))S(v)d,
\]
\begin{multline*}
S(AS(B)cS(d))e=AS(B)cS(S(d)e)=AS(B)S(cS(d))e\\
=AS(S(B)cS(d))e=S(AS(B))cS(d)e.
\end{multline*}
\end{Proof}

\begin{lemma}
\label{lem:T-L}
For every $k\geq4$, $\LL(\phi(T_k))$ and
$\RR(\phi(\widetilde T_k))$ hold.
\end{lemma}
\begin{Proof}
For $k=4$, we have
$$\phi(T_4)=S(x_1)S(S(x_2)x_3)x_4
\sim_S S(S(x_1)S(x_2))x_3x_4\sim_S S(x_1)S(S(x_2)x_3x_4).$$
If $u=AS(B)$, then $S(u)z\sim_S AS(S(B)z)$. Thus $\LL(u)$ implies $\LL(S(u)z)$, and induction proves the assertion for $T_k$. The other assertion follows by reflection.
\end{Proof}

\begin{theorem}
\label{thm:validity}
The identities $F_r,\rho F_r,G_r,\rho G_r,J_r,\rho J_r,H_{p,q}$ hold
for all the stated parameters.
\end{theorem}
\begin{Proof}
For $F_r$, put $P_1=c_1$, $P_j=S(P_{j-1})c_j$ for $j\geq2$, and
$W_2=S(b)S(c_1)c_2$, $W_j=S(W_{j-1})c_j$ for $j\geq3$. We prove
\begin{equation}
S(a)S(b)S(P_j)=S(a)S(W_j),\qquad j\geq2.
\label{Finduction}
\end{equation}
For $j=2$, relation \eqref{R} gives
$$S(a)S(b)S(S(c_1)c_2)=S(S(a)S(b))S(c_1)c_2=S(a)S(S(b)S(c_1)c_2)=S(a)S(W_2).$$
Suppose \eqref{Finduction} holds for $j$. Then
\begin{multline*}
S(a)S(b)S(P_{j+1})
=S(a)S(b)S(S(P_j)c_{j+1})
=S(S(a)S(b)S(P_j))c_{j+1}\\
=S(S(a)S(W_j))c_{j+1}
=S(a)S(S(W_j)c_{j+1})=S(a)S(W_{j+1}).
\end{multline*}
Taking $j=r-1$ and multiplying by $c_r$ gives the operated images of $F_r$.

For $G_r$, let $B_r=\phi(L_r(b_1,\ldots,b_r))$,
$A_r=\phi(T_{r+1}(a,b_1,\ldots,b_r))$, and $U_r=S(a)B_r$.
For $r\geq3$, we have
\[
U_r\sim_S S(S(a)S(B_{r-2}))b_{r-1}b_r
\sim_S S(a)S(S(B_{r-2})b_{r-1}b_r),
\]
so $\LL(U_r)$ holds. By \eqref{J0} and \eqref{H0},
\begin{equation}\label{Gaux}
S(a)S(B_r)S(c)d=S(S(a)S(B_{r-1}))b_rS(c)d=S(S(a)S(B_{r-1})b_rS(c))d=S(U_r)S(c)d.
\end{equation}
Since $A_3=U_3$, this proves $G_3$. For $r\geq4$, use
$A_r=S(A_{r-1})b_r$ and the induction hypothesis to obtain
\begin{multline*}
S(A_r)S(c)d=S(A_rS(c))d
=S(A_{r-1})b_rS(S(c)d)
=S(A_{r-1})S(b_rS(c))d\\
=S(a)S(B_{r-1})S(b_rS(c))d
=S(a)S(B_{r-1})b_rS(S(c)d)\\
=S(U_rS(c))d=S(U_r)S(c)d=S(a)S(B_r)S(c)d.
\end{multline*}
The first and penultimate equalities use \eqref{H0}; the last uses
\eqref{Gaux}. Hence $G_r$ holds.

To prove $H_{p,q}$, we first establish
\begin{equation}
uS(zS(v))=S(u)zS(v),\qquad \LL(u),\ \RR(v).
\label{generalB}
\end{equation}
Choose $v=S(C)D$. By \eqref{R} and \eqref{H0},
\begin{multline*}
uS(zS(v))=uS(zS(S(C)D))=uS(S(zS(C))D)\\
=S(uS(zS(C)))D=S(u)S(zS(C))D
=S(u)zS(S(C)D)=S(u)zS(v).
\end{multline*}
Now substitute $u=\phi(T_p)$ and $v=\phi(\widetilde T_q)$ and use
Lemma~\ref{lem:T-L}.

For $J_r$, suppose $\RR(v)$ and choose $v=S(P)Q$. We claim that
\begin{equation}
A S(b)cS(dS(v))=A S(S(b)c)dS(v).
\label{generalE}
\end{equation}
Indeed,
\begin{multline*}
A S(b)cS(dS(v))=A S(b)cS(dS(S(P)Q))
=A S(b)cS(S(dS(P))Q)\\
=A S(b)S(cS(dS(P)))Q
=A S(S(b)cS(dS(P)))Q
=S(A S(b))cS(dS(P))Q,
\end{multline*}
whereas
$$A S(S(b)c)dS(v)=A S(S(b)c)dS(S(P)Q)=A S(S(b)c)S(dS(P))Q=S(A S(b))cS(dS(P))Q.$$
Taking $A=S(a)$ and $v=\phi(\widetilde T_r)$ proves $J_r$. Reflection preserves \eqref{R}, so it also proves the reflected identities.
\end{Proof}

\noindent\textbf{Regular operated words.}

Let $X$ be the set of variables.

\begin{definition}
A variable is a regular word. An operated word is called \emph{regular}
if it has the form
\[
S(u_1)\cdots S(u_p)xS(v_1)\cdots S(v_q),
\]
where $x\in X$, $p,q\geq 0$, and all $u_i,v_j$ are regular.
\end{definition}

\begin{lemma}
\label{lem:regular}
An operated word is regular if and only if it is the image under $\phi$
of a monomial in $\prec,\succ$.
\end{lemma}

\begin{Proof}
Since $\phi(U\prec V)=\phi(U)S(\phi(V))$ and
$\phi(U\succ V)=S(\phi(U))\phi(V)$, the image of every monomial is
regular.

Conversely, let $w=S(u_1)\cdots S(u_p)xS(v_1)\cdots S(v_q)$ be regular.
By induction on the number of variables, we may choose monomials
$U_i,V_j$ such that $\phi(U_i)=u_i$ and $\phi(V_j)=v_j$. Then $w$ is
the image of
\[
\bigl(\cdots\bigl(
U_1\succ(U_2\succ\cdots\succ(U_p\succ x))
\bigr)\prec V_1\cdots\bigr)\prec V_q,
\]
where the corresponding products are omitted if $p=0$ or $q=0$.
\end{Proof}

\begin{lemma}
\label{lem:literal-unique}
If $\phi(U)=\phi(V)$ as operated words, then $U=V$ follows from
\[
(A\succ B)\prec C=A\succ(B\prec C).
\]
Moreover, for every occurrence of a regular subword $u$ in a regular
word $w$, there exists a monomial $W$ with $\phi(W)=w$ in which this
occurrence of $u$ is the image of a submonomial.
\end{lemma}

\begin{Proof}
Proceed by induction on the number of variables in $w=\phi(U)=\phi(V)$.
Use the decomposition of $w$ from Lemma~\ref{lem:regular}. If a preimage
of $w$ has the form $A\prec B$, then $\phi(B)=v_q$ and
$\phi(A)=S(u_1)\cdots S(u_p)\allowbreak xS(v_1)\cdots S(v_{q-1})$.
If it has the form $A\succ B$, then $\phi(A)=u_1$ and
$\phi(B)=S(u_2)\cdots S(u_p)\allowbreak xS(v_1)\cdots S(v_q)$.
The given identity allows us to interchange adjoining a factor on
the left by $\succ$ and adjoining a factor on the right by $\prec$.
Thus, by induction, every preimage is equal to the monomial constructed
in Lemma~\ref{lem:regular}.

Suppose first that the chosen occurrence of $u$ is not contained in
an argument of $S$. Then $u$ contains $x$ and consists of consecutive
factors of $w$. Choose a preimage of $u$ and adjoin the remaining
factors on the left by $\succ$ and on the right by $\prec$.
If $u$ lies in an argument of $S$, apply the same construction
inductively to that argument. Hence the chosen occurrence of $u$
is the image of a submonomial. The first assertion shows that
replacing this submonomial is independent, modulo the given identity,
of the preimages chosen.
\end{Proof}

\begin{lemma}
\label{lem:hiddenarm}
Let $T=A\succ B$ be a monomial containing only $\succ$. Then
\[
\LL(\phi(T))
\quad\Longleftrightarrow\quad
\LL(\phi(A))\ \text{or}\ B\text{ is not a right }\succ\text{-comb},
\]
where a variable is also regarded as a right $\succ$-comb.
If these equivalent conditions hold, then $T$ is equal modulo $F_3$
to $T_k(U_1,\ldots,U_k)$ for some $k\geq4$, where all $U_i$
contain only $\succ$. Otherwise, no application of \eqref{R}
in either direction is possible in $\phi(T)$.
The corresponding statements for $\prec$ and $Q_k$ follow by reflection.
\end{lemma}

\begin{Proof}
Proceed by induction on the number of variables.
If $\LL(\phi(A))$ holds, induction gives
$A=T_k(U_1,\ldots,U_k)$ modulo $F_3$. Hence
$T=T_{k+1}(U_1,\ldots,U_k,B)$ modulo $F_3$.

Suppose that $B$ is not a right $\succ$-comb. Then
\[
T=A\succ(D_1\succ\cdots\succ(D_m\succ((C\succ D)\succ E)))
\]
for some $m\geq0$. If $m=0$, then $T=T_4(A,C,D,E)$.
For $m>0$, apply $F_3$ in the form
\[
P\succ\bigl(Q\succ((C\succ D)\succ E)\bigr)
=
P\succ\bigl((Q\succ(C\succ D))\succ E\bigr).
\]
Applied to the innermost such submonomial, this decreases $m$ by one.
Repeating this transformation gives a substitution into $T_4$.
In both cases, Lemma~\ref{lem:T-L} implies $\LL(\phi(T))$,
since its conclusion is preserved under substitution.

Now suppose that $\LL(\phi(A))$ does not hold and that $B$ is a
right $\succ$-comb. Write
$B=x_1\succ(x_2\succ\cdots\succ x_m)$, allowing $m=1$. Then
\[
\phi(T)=S(\phi(A))S(x_1)\cdots S(x_{m-1})x_m.
\]
By induction, no application of \eqref{R} is possible inside
$\phi(A)$. At the outer level, a subword of the form $US(S(V)W)$
cannot occur: the only possible nonsingleton $S$-argument is
$\phi(A)$, and its $S$-factor is the first factor of the word.
A subword of the form $S(US(V))W$ cannot occur either, since
$\phi(A)$ ends in a variable. Thus no application of \eqref{R}
in either direction is possible. Since $\phi(T)$ ends in $x_m$,
its class contains no word of the form $US(V)$ with $U,V$ nonempty.
Therefore $\LL(\phi(T))$ does not hold.

Reflection proves the assertions for $\prec$.
\end{Proof}

\begin{lemma}
\label{lem:rootfactor}
Every monomial containing $\prec$ is equal modulo \eqref{NB}
to a monomial of the form $A\prec B$.
Similarly, every monomial containing $\succ$ is equal modulo
\eqref{NB} to a monomial of the form $A\succ B$.
\end{lemma}

\begin{Proof}
Proceed by induction on the number of variables.
It suffices to consider a monomial $U\succ V$ containing $\prec$.
If $U$ contains $\prec$, write $U=A\prec B$ modulo \eqref{NB}
by induction. Then
$(A\prec B)\succ V=A\prec(B\succ V)$.
Otherwise, write $V=A\prec B$ and use
$U\succ(A\prec B)=(U\succ A)\prec B$.
The second assertion follows by interchanging $\prec$ and $\succ$.
\end{Proof}
Combining Lemmas~\ref{lem:hiddenarm} and \ref{lem:rootfactor}, every binary
monomial $U$ with $\LL(\phi(U))$ satisfies, modulo \eqref{NB} and $F_3$,
\begin{equation*}
 U\equiv A\prec B
 \quad\text{or}\quad
 U\equiv T_p(U_1,\ldots,U_p),\ p\ge4,
\end{equation*}
and, when $\RR(\phi(V))$ holds, the reflected identities give
\begin{equation*}
 V\equiv C\succ D
 \quad\text{or}\quad
 V\equiv \widetilde T_q(V_1,\ldots,V_q),\ q\ge4.
\end{equation*}

\noindent\textbf{Singleton intervals.}

Let $w$ be a multilinear operated word in which the variables
$x_1,\ldots,x_n$ occur in this order. To each occurrence of $S$
we associate the interval $[i,j]$ of variables in its argument.
Denote the resulting family of intervals by $\Ical(w)$, counting
each occurrence separately. Any two intervals in $\Ical(w)$
are disjoint or one contains the other.

Suppose that the factors $A,B,C$ in \eqref{R} have supports
$[\ell,p-1]$, $[p,q]$, and $[q+1,r]$, respectively. Then
\eqref{R} replaces the interval $[p,r]$ by $[\ell,q]$, or conversely,
while preserving $[p,q]$:
\begin{equation}
\{[p,q],[p,r]\}\longleftrightarrow
\{[p,q],[\ell,q]\},
\qquad \ell<p\leq q<r.
\label{intervalmove}
\end{equation}
All other intervals remain unchanged. We call $[p,q]$ the pivot
of this transformation.

Put
\begin{equation*}
M(w)=\{i:[i,i]\in\Ical(w)\}.
\end{equation*}

\begin{lemma}
\label{lem:singleton}
If $w\sim_S w'$, then $M(w)=M(w')$. Moreover, if $w$ is a derived word, then every interval in $\Ical(w')$ contains an element of $M(w)$.
\end{lemma}

\begin{Proof}
Both intervals replaced in \eqref{intervalmove} properly contain the pivot. Hence neither is a singleton, and $M(w)$ is preserved by every application of \eqref{R}.

In a derived word, every minimal occurrence of $S$ has the form $S(x_i)$. Thus every interval in $\Ical(w)$ contains an element of $M(w)$. Under \eqref{intervalmove}, the new interval contains the unchanged pivot. Consequently, this property is preserved by each application of \eqref{R}, proving the second assertion.
\end{Proof}

For a binary monomial $T$, write $M(T)=M(\phi(T))$ and put
\[
\Ncal_n(M)=\{T\in\Ncal_n:M(T)=M\}.
\]
We call $\Ncal_n(M)$ the sector with singleton set $M$.
The singleton invariant gives a decomposition of the normal basis
according to the subsets $M\subseteq\{1,\ldots,n\}$. Define
\[
K_n(M)=K_n\cap\Span\Ncal_n(M),
\qquad
Q_n(M)=K_n(M)/(D_n\cap K_n(M)).
\]
Then
\begin{equation}
K_n=\bigoplus_M K_n(M),
\qquad
Q_n=\bigoplus_M Q_n(M).
\label{sector-decomp}
\end{equation}

Indeed, by Lemma~\ref{lem:fibres}, $K_n$ is spanned by differences
of normal monomials with $S$-equivalent images. Each such difference
belongs to a single $K_n(M)$, which proves the first decomposition.

The space $D_n$ is spanned by consequences of binomial identities
of smaller arity. Operadic substitutions and multiplication by
binary monomials preserve the $S$-equivalence of their images.
Reduction to normal monomials modulo \eqref{NB} also preserves
this equivalence. Hence $D_n$ is spanned by binomials, each
belonging to one $K_n(M)$. Therefore
$D_n=\bigoplus_M(D_n\cap K_n(M))$, and the second decomposition
follows.

\begin{definition}
Let $T$ be a binary monomial in the ordered variables
$x_1,\ldots,x_n$. Define its \emph{center} recursively by
\[
c(x_i)=i,\qquad
c(U\prec V)=c(U),\qquad
c(U\succ V)=c(V).
\]
For each $i$, put
\begin{equation*}
d_i(T)=
\#\{W:\ W\text{ is a submonomial of }T,\ \deg W\geq2,\ c(W)=i\},
\end{equation*}
where submonomials are counted by their occurrences.
We call $d(T)=(d_1(T),\ldots,d_n(T))$ the \emph{degree vector}
of $T$.
\end{definition}

Since a binary monomial of degree $n$ has exactly $n-1$ internal vertices, we have
\begin{equation*}
d_i(T)\geq0,\qquad
\sum_{i=1}^n d_i(T)=n-1.
\end{equation*}

\begin{remark}
\label{rem:degree-not-invariant}
The degree vector need not be preserved by \eqref{NB} in a larger context.
For example,
\[
x_1\succ\bigl(x_2\prec(x_3\succ x_4)\bigr)
=x_1\succ\bigl((x_2\prec x_3)\succ x_4\bigr),
\]
but the two degree vectors are $(0,2,0,1)$ and $(0,1,0,2)$.
A defining relation can change the center of a submonomial and hence the
centers of the monomials containing it. We shall use degree vectors as
coordinates on the fixed normal basis. Their zero sets, however, are
invariant, as Lemma~\ref{lem:center-tree} will show.
\end{remark}

\begin{lemma}
\label{lem:pure-vectors}
Let $a=(a_1,\ldots,a_m)$ be a nonnegative integer vector with
$\sum_i a_i=m-1$. It is the degree vector of a unique pure $\prec$-monomial
if and only if
\begin{equation}
\sum_{i=1}^t a_i\geq t,\qquad 1\leq t<m.
\label{prec-vector}
\end{equation}
It is the degree vector of a unique pure $\succ$-monomial if and only if
\begin{equation*}
\sum_{i=1}^t a_i<t,\qquad 1\leq t<m.
\end{equation*}
\end{lemma}
\begin{Proof}
In a pure $\prec$-monomial, the center of each submonomial is its first
leaf. We reconstruct the monomial from right to left. At position $i$,
start with $x_i$ and attach the nearest $a_i$ previously constructed
monomials, successively, as right factors of $\prec$. Replace these
monomials by the resulting one.

Put $D_i=\sum_{j\leq i}a_j$. Before position $i$ is processed, the
number of available monomials is
$(m-i)-\sum_{j>i}a_j=D_i-i+1$. This is at least $a_i$ precisely when
$D_{i-1}\geq i-1$. Thus \eqref{prec-vector} makes the construction
possible, and the total sum leaves one monomial at the end. Exactly
$a_i$ vertices with center $i$ have been added. Conversely, removing the
successive roots with center $i$ from a pure monomial recovers this
construction. This proves necessity and uniqueness.

Reflection reverses the vector. Since its total sum is $m-1$, the
reflected inequalities are $D_t\leq t-1$ for $t<m$, which proves the
second assertion.
\end{Proof}

For a weak composition $d=(d_1,\ldots,d_n)$ of $n-1$, put
$D_j=\sum_{i=1}^j d_i$ and $h_j=D_j-j$, with $D_0=h_0=0$.

\begin{theorem}
\label{thm:degreevector}
The map $T\mapsto d(T)$ is a bijection from $\Ncal_n$ to
\[
\{(d_1,\ldots,d_n)\in\mathbb Z_{\geq0}^n:d_1+\cdots+d_n=n-1\}.
\]
Moreover,
\begin{equation}
c(T)=1+\max\{j<n:D_j<j\},
\label{centerformula}
\end{equation}
where the maximum is taken to be $0$ if the set is empty.
\end{theorem}
\begin{Proof}
Proceed by induction on $n$. There is nothing to prove for $n=1$.

Suppose first that $d_n>0$. A normal monomial $U\prec P$ has degree
zero at its last leaf, so the required root must be $U\succ Q$, with
$Q$ pure $\succ$. Put
\begin{equation}
k=\max\{j<n:D_j=j-1\}.
\label{k-succ}
\end{equation}
This set is nonempty: $h_{n-1}=-d_n<0$, whereas $h_0=0$ and
$h_j-h_{j-1}=d_j-1\geq-1$. By the maximality of $k$,
\begin{equation*}
h_j\leq-2,\qquad k<j<n.
\end{equation*}
Indeed, a later value at least $-1$ would force another occurrence of
$-1$ before the final negative value.

The prefix $(d_1,\ldots,d_k)$ has sum $k-1$ and determines $U$ by
induction. For the suffix, put
\[d^Q=(d_{k+1},\ldots,d_{n-1},d_n-1).\]
Its sum is $n-k-1$. Since $h_k=-1$, its proper prefix sums satisfy
\[D_{k+t}-D_k=t+h_{k+t}-h_k<t.\]
Lemma~\ref{lem:pure-vectors} therefore gives a unique pure monomial $Q$ with vector $d^Q$.
The root of $U\succ Q$ adds one at the last leaf, giving the vector $d$.
Conversely, the prefix inequalities for a pure right factor force the
split \eqref{k-succ}. Hence the reconstruction is unique.

Now suppose $d_n=0$. Then $h_{n-1}=0$. Put
\begin{equation*}
r=1+\max\{j<n:h_j<0\},\qquad
k=\min\{j\geq r:h_j=0\},
\end{equation*}
with $r=1$ if the first set is empty. We have
\begin{equation}
h_j>0\quad(r\leq j<k),\qquad
h_j\geq0\quad(k\leq j<n).
\label{prec-h}
\end{equation}
In particular, $d_r>0$. Subtract one at position $r$ in the prefix
$(d_1,\ldots,d_k)$. The resulting vector has sum $k-1$ and determines
a unique normal monomial $U$. By the induction hypothesis and
\eqref{prec-h}, its center is $r$.

The suffix $(d_{k+1},\ldots,d_n)$ has sum $n-k-1$, and each proper
prefix sum is $D_{k+t}-D_k=t+h_{k+t}\geq t$. It therefore determines
a unique pure $\prec$-monomial $P$. The root of $U\prec P$ adds one
at $c(U)=r$, restoring $d$. Conversely, the center formula for $U$ and
the prefix inequalities for $P$ force these choices of $r$ and $k$.
Thus this reconstruction is also unique.

In the first case, $c(T)=n$ and $D_{n-1}<n-1$. In the second,
$c(T)=r$ by construction. This proves \eqref{centerformula} and the
bijection.
\end{Proof}

In particular, $\dim\;\NB(n)=\binom{2n-2}{n-1}$.

\noindent\textbf{The center tree and singleton positions.}

Associate a directed tree $\Gamma(T)$ with a binary monomial $T$.
At a root $U\prec V$, add the edge $c(U)\to c(V)$; at a root
$U\succ V$, add $c(V)\to c(U)$. Retain the edges of the two subtrees.
The root of $\Gamma(T)$ is $c(T)$.

\begin{lemma}
\label{lem:center-tree}
Let $T$ have arity $n\geq2$. Then
$d_i(T)=\outdeg_{\Gamma(T)}(i)$. Each occurrence of $S$ in $\phi(T)$
corresponds to an edge of $\Gamma(T)$, and its argument is supported on
the descendant subtree rooted at the target of that edge. Consequently,
\begin{equation}
M(T)=\{i:d_i(T)=0\}.
\label{M-zero-set}
\end{equation}
\end{lemma}
\begin{Proof}
Each internal vertex of $T$ contributes one edge starting at its center.
This gives the outdegree formula. For $U\prec V$, the new operator is
applied to $\phi(V)$, whose leaves form the descendant subtree rooted at
$c(V)$. For $U\succ V$, the same argument applies to $U$. Induction
proves the correspondence for all occurrences of $S$.

A singleton interval occurs precisely at a leaf of $\Gamma(T)$. Since
$n\geq2$, such a vertex has an incoming edge and no outgoing edges.
This proves \eqref{M-zero-set}.
\end{Proof}

\begin{corollary}
\label{cor:sector-size}
If $n\geq2$ and $M\subseteq\{1,\ldots,n\}$ has cardinality $s$,
where $1\leq s\leq n-1$, then
\begin{equation*}
|\Ncal_n(M)|=\binom{n-2}{s-1}.
\end{equation*}
\end{corollary}
\begin{Proof}
By \eqref{M-zero-set}, put $d_i=0$ for $i\in M$ and $d_i=1+e_i$
otherwise. Then $\sum_{i\notin M}e_i=s-1$. There are $n-s$ nonnegative
coordinates, so the number of these vectors is $\binom{n-2}{s-1}$.
Apply Theorem~\ref{thm:degreevector}.
\end{Proof}

At a nonsingleton position $i$, we call $d_i-1$ its excess.
When listing excess positions, we count position $i$ with multiplicity
$d_i-1$.

\section{The corner and the pure case}
\label{sec:corner}

For $m\geq3$ and $2\leq h\leq m-1$, put
\begin{equation}
P_m^h=L_{m-h}\bigl(x_1\succ L_h(x_2,\ldots,x_{h+1}),
                    x_{h+2},\ldots,x_m\bigr).
\label{Pmh}
\end{equation}
Let $R_j$ denote the right $\prec$-comb on the indicated consecutive
leaves. When the argument of $\Ical$ is a binary monomial, we mean its
operated image. Thus
\begin{equation*}
\Ical(P_m^h)=\{[1,1]\}\cup\{[2,j]:2\leq j\leq h\}
\cup\{[1,j]:h+1\leq j\leq m-1\}.
\end{equation*}

\begin{lemma}\label{lem:pair-barrier}
No word $S$-equivalent to $\phi(P_m^h)$ contains the interval $[1,h]$.
The assertion remains true after adjoining fixed intervals strictly to
the right of this prefix, provided every possible pivot contains
position $2$ and no other interior singleton.
\end{lemma}
\begin{Proof}
In the sector $\{1,2\}$, a pivot cannot contain the first position and
must contain a singleton. Hence it starts at $2$. Creating $[1,h]$
requires the intervals $[2,h]$ and $[2,r]$ with $r>h$. They do not
occur together initially. Suppose they first occur together before
$[1,h]$ is created.

If the new interval is $[2,r]$, it must arise from
\[
\{[2,q],[1,q]\}\longrightarrow\{[2,q],[2,r]\},\qquad q<r,
\]
while $[2,h]$ remains unchanged. For $q<h$, the old interval $[1,q]$
crosses $[2,h]$. For $q=h$, the interval $[1,h]$ is already present.
For $q>h$, the intervals $[2,h]$ and $[2,q]$ already form the required
pair. Each case is impossible.

If $[2,h]$ is new and $[2,r]$ is unchanged, the move requires an old
interval $[1,q]$ with $q<h$, which crosses $[2,r]$. Thus the pair
cannot occur, and $[1,h]$ cannot be created. Fixed intervals to the
right do not change this argument under the stated pivot condition.
\end{Proof}

\begin{corollary}
\label{cor:pair-zero}
For $m\geq3$, $K_m(\{1,2\})=K_m(\{m-1,m\})=0$.
\end{corollary}
\begin{Proof}
The vector of $P_m^h$ has zeros at $1,2$, one at every other position,
and one additional unit at $h+1$. These are all vectors in the first
sector. If $h'<h$, then $P_m^{h'}$ contains $[1,h]$, whereas no word
congruent to $P_m^h$ does. Hence all these normal monomials lie in
distinct fibres. Reflection gives the other sector.
\end{Proof}

\begin{proposition}\label{cor:triple-zero-kernel}
For $n\geq4$, $K_n(\{1,2,n\})=K_n(\{1,n-1,n\})=0$.
\end{proposition}
\begin{Proof}
The normal monomials in the first sector are
\begin{equation*}
P_r^h\prec R_{n-r}(x_{r+1},\ldots,x_n),
\qquad 3\leq r\leq n-1,\quad 2\leq h\leq r-1.
\end{equation*}
Their two excess positions are $h+1,r$, which range over all unordered
pairs, with repetition, from $\{3,\ldots,n-1\}$.
Every pivot avoids the two endpoints and contains a singleton, so it
starts at $2$. All intervals not containing $2$ are therefore fixed.
In particular, the intervals $[r+1,n],\ldots,[n,n]$ determine $r$ and
remain unchanged along every path. For a fixed $r$,
Lemma~\ref{lem:pair-barrier} separates the different values of $h$.
Thus each fibre contains one normal monomial. Apply reflection for the
second sector.
\end{Proof}

\begin{lemma}\label{lem:rigid-boundary}
Let $u$ be regular, and suppose no application of \eqref{R} in either
direction is possible inside $u$. If $u$ ends in a direct variable and
$S(u)$ is a block on a global initial interval, this block remains
unchanged along every $S$-path. Replacing its argument $u$ by one new
variable gives a valid path. The reflected assertion holds for a final
block $S(u)$ when $u$ starts in a direct variable.
\end{lemma}
\begin{Proof}
Consider the first move changing the block. It cannot be internal.
Removing $S(u)$ as the right interval in \eqref{intervalmove} requires
a nonempty factor before the first leaf. Removing it as the left
interval requires $u$ to end in an $S$-factor. Both are impossible.
If the block is the pivot, its argument remains unchanged. Hence every
move respects the block and descends after replacing $u$ by a variable.
The endpoints remain regular, and Lemma~\ref{lem:literal-unique}
identifies expansion with binary substitution.
\end{Proof}

We call such a block rigid. In every contraction below, we retain its
outer operator: $S(u)$ is replaced by $S(y)$. Thus a block of $t$ leaves
loses $t-1$ leaves and $t-1$ operators, and the endpoints remain regular.

Put $W_3=x_1\succ(x_2\succ x_3)$ and $W_{j+1}=W_j\succ x_{j+1}$.
The words $\phi(W_j)$ are internally irreducible in both directions.
Indeed, $\phi(W_3)=S(x_1)S(x_2)x_3$, and
$\phi(W_{j+1})=S(\phi(W_j))x_{j+1}$. Every nontrivial $S$-factor is
first on its level and its argument ends in a direct variable.

\noindent\textbf{The four-singleton corner.}

Let $n\geq7$, put $M_0=\{1,2,n-1,n\}$ and $K=n-4$, and define
\begin{equation*}
\Delta_K=\{(a,b,c):1\leq a\leq b<c\leq K+1\}.
\end{equation*}
For $(a,b,c)\in\Delta_K$, let $r=b+2$ and $h=r-a$. Define
\begin{equation*}
\Theta(a,b,c)=
\begin{cases}
(P_r^h\prec R_{n-r-1})\prec x_n,&c=b+1,\\
P_r^h\prec\rho(P_{n-r}^{\,n-c-1}),&c\geq b+2,
\end{cases}
\end{equation*}
where all factors use consecutive leaves. These monomials are normal.
Their three excess positions, counted with multiplicity, are
\begin{equation}
u=b-a+3,\qquad v=b+2,\qquad w=c+1.
\label{corner-excess}
\end{equation}
This follows from \eqref{Pmh}, since each root $\prec$ adds one at the
last leaf of $P_r^h$. Conversely, $b=v-2$, $a=v-u+1$ and $c=w-1$.
Thus $\Delta_K$ parametrizes $\Ncal_n(M_0)$.

\begin{lemma}
\label{lem:corner-locked}
The fibre of $\Theta(a,b,c)$ contains no other normal monomial if
\begin{equation}
a=b\geq2\qquad\text{or}\qquad c=K+1,\ b\leq K-1.
\label{corner-locked-shapes}
\end{equation}
\end{lemma}
\begin{Proof}
In the first case, $h=2$ and $P_r^2=W_r$, with $r\geq4$.
Its image begins with $S(\phi(W_{r-1}))$. By
Lemma~\ref{lem:rigid-boundary}, we may contract this argument along
any path to a regular endpoint. The contracted endpoints have arity
$m=n-r+2$ and singleton set $\{1,m-1,m\}$. Their binary preimages
are equal by Proposition~\ref{cor:triple-zero-kernel}. Expanding the
same argument and applying Lemma~\ref{lem:literal-unique} proves
equality of the original normal monomials.

In the second case, the right factor is $\rho(P_q^2)=\rho(W_q)$,
where $q=n-r\geq3$. Its image is the argument of a rigid final
$S$-block. Contracting it gives the sector $\{1,2,m\}$ with $m=r+1$.
The same argument applies.
\end{Proof}

We call the coordinates in \eqref{corner-locked-shapes} locked. Their
complement is
\begin{equation}
\begin{split}
\mathcal M_K={}&\{(1,1,c):2\leq c\leq K\}\\
&\cup\{(a,b,c):1\leq a<b<c\leq K\}\\
&\cup\{(a,K,K+1):1\leq a\leq K-1\}.
\end{split}
\label{MK}
\end{equation}

\begin{lemma}
\label{lem:corner-embeddings}
Inserting $\succ$ at leaf $2$ or $\prec$ at leaf $n-2$ of an
arity-$(n-1)$ corner monomial preserves normality. The coordinate maps
are
\begin{equation}
\iota_L(a,b,c)=(a,b+1,c+1),\qquad
\iota_R(a,b,c)=(a,b,c).
\label{corner-embeddings}
\end{equation}
\end{lemma}
\begin{Proof}
Put $m=n-1$. A nontrivial pure $\prec$-factor containing leaf $2$
would start at $1$ or $2$. Its first leaf would have positive degree,
contrary to $d_1=d_2=0$. Thus inserting $\succ$ at leaf $2$ preserves
the normal grammar. Similarly, a nontrivial pure $\succ$-factor
containing leaf $m-1$ would end at $m-1$ or $m$, contrary to the zero
degrees there. This proves normality for the second insertion.
The first insertion increases all three excess positions by one; the
second leaves them unchanged. Formula~\eqref{corner-excess} gives
the stated maps.
\end{Proof}

For $n=7$, the coordinates of $\mathcal M_3$ are
$(1,1,2)$, $(1,1,3)$, $(1,2,3)$, $(1,3,4)$ and $(2,3,4)$.
The following four paths connect them. The letters $a,b,c,d,e,f,g$
denote the leaves in order.

For $(1,1,2)$ and $(1,2,3)$, we have
\begin{multline*}
S(a)S(b)cS(dS(eS(f)))S(g)
\sim_S S(a)S(b)cdS(S(eS(f))S(g))\\
\sim_S S(a)S(S(b)cdS(eS(f)))S(g)
\sim_S S(S(a)S(b))cdS(eS(f))S(g)\\
\sim_S S(a)S(S(b)c)dS(eS(f))S(g).
\end{multline*}
For $(1,2,3)$ and $(1,3,4)$,
\begin{multline*}
S(a)S(S(b)c)dS(eS(f))S(g)
\sim_S S(a)S(S(b)c)deS(S(f)S(g))\\
\sim_S S(a)S(S(S(b)c)deS(f))S(g)
\sim_S S(S(a)S(S(b)c))deS(f)S(g)\\
\sim_S S(a)S(S(S(b)c)d)eS(f)S(g).
\end{multline*}
For $(1,2,3)$ and $(2,3,4)$,
\begin{multline*}
S(a)S(S(b)c)dS(eS(f))S(g)
\sim_S S(S(a)S(b))cdS(eS(f))S(g)\\
\sim_S S(a)S(S(b)cd)S(eS(f))S(g)
\sim_S S(a)S(S(b)cd)eS(S(f)S(g))\\
\sim_S S(a)S(S(S(b)cd)eS(f))S(g)
\sim_S S(S(a)S(S(b)cd))eS(f)S(g)\\
\sim_S S(S(S(a)S(b))cd)eS(f)S(g)
\sim_S S(S(a)S(S(b)c)d)eS(f)S(g).
\end{multline*}
For $(1,2,3)$ and $(1,1,3)$,
\begin{multline*}
S(a)S(S(b)c)dS(eS(f))S(g)
\sim_S S(S(a)S(b))cdS(eS(f))S(g)\\
\sim_S S(S(a)S(b))cdeS(S(f)S(g))
\sim_S S(S(a)S(b))cS(deS(f))S(g)\\
\sim_S S(a)S(S(b)cS(deS(f)))S(g)
\sim_S S(a)S(b)cS(S(deS(f))S(g))\\
\sim_S S(a)S(b)cS(deS(S(f)S(g)))
\sim_S S(a)S(b)cS(dS(eS(f))S(g)).
\end{multline*}

\begin{proposition}
\label{prop:corner-mobile}
For $n\geq7$, the coordinates in $\mathcal M_K$ form one $S$-fibre.
Every other fibre in the corner contains one normal monomial.
\end{proposition}
\begin{Proof}
For $n=7$, the four paths above prove the first assertion, and
Lemma~\ref{lem:corner-locked} proves the second.

Suppose $K\geq4$ and the assertion holds for $K-1$. The set
\begin{equation*}
\mathcal O_K=\iota_R(\mathcal M_{K-1})\cup
\iota_L(\mathcal M_{K-1})
\end{equation*}
is connected. Indeed, each image is connected, and both contain
$(1,2,3)$, with preimages $(1,2,3)$ and $(1,1,2)$, respectively.
From \eqref{MK}, its complement in $\mathcal M_K$ is
\begin{equation*}
\mathcal E_K=\{(1,1,K)\}\cup
\{(t,t+1,K):2\leq t\leq K-3\}\cup\{(K-1,K,K+1)\}.
\end{equation*}
The special identities give the edges
\begin{align}
J_{n-4}:&\quad (1,1,K)\longleftrightarrow(1,2,K),\label{Jcorner}\\
H_{p,n-p-1}:&\quad(p-3,p-2,K)\longleftrightarrow(p-2,p-1,K),
\quad4\leq p\leq K-1,\notag\\
\rho J_{n-4}:&\quad(K-1,K,K+1)\longleftrightarrow(K-2,K-1,K).
\label{rhoJcorner}
\end{align}
For the first two, the excess positions are
$(3,3,n-3)\leftrightarrow(4,4,n-3)$ and
$(4,p,n-3)\leftrightarrow(4,p+1,n-3)$. For the last, normalize the
reflected endpoints by \eqref{NB} to obtain
\[
(T_{n-2}\prec x_{n-1})\prec x_n,
\qquad (T_{n-3}\prec(x_{n-2}\prec x_{n-1}))\prec x_n.
\]
Their excess positions give \eqref{rhoJcorner}.

The first and last edges join the two extreme points of $\mathcal E_K$
to $\mathcal O_K$. Starting at $(1,2,K)\in\mathcal O_K$, the $H$-edges
join the remaining points successively. Thus $\mathcal M_K$ is one
fibre. Lemma~\ref{lem:corner-locked} excludes all other coordinates.
\end{Proof}

\begin{corollary}
\label{cor:corner-upper}
For $n\geq8$, the corner kernel modulo $D_n$ is spanned by
$J_{n-4},\rho J_{n-4}$ and $H_{p,n-p-1}$, $4\leq p\leq n-5$.
\end{corollary}
\begin{Proof}
Differences within either image of $\mathcal M_{K-1}$ are consequences
of identities of smaller arity. Since the images intersect, all
differences within $\mathcal O_K$ belong to $D_n$. The edges
\eqref{Jcorner}--\eqref{rhoJcorner} join every remaining point of
$\mathcal M_K$ to this set. The other fibres contribute no kernel.
\end{Proof}

\noindent\textbf{The pure $\succ$ case.}
Let $\Mag$ be the free nonsymmetric operad generated by $\succ$ and let
\[
 \psi:\Mag\longrightarrow\AsS,\qquad
 \psi(U\succ V)=S(\psi(U))\psi(V).
\]
Besides $F_r$ and $G_r$, we use two identities of arity six:
\begin{equation*}
\begin{split}
K_6:\quad
\bigl(a\succ(b\succ((c\succ d)\succ e))\bigr)\succ f
&=a\succ\bigl(b\succ(((c\succ d)\succ e)\succ f)\bigr),\\
E_6:\quad
\bigl(a\succ((b\succ c)\succ(d\succ e))\bigr)\succ f
&=a\succ\bigl(((b\succ c)\succ(d\succ e))\succ f\bigr).
\end{split}
\end{equation*}
The operated images of the two sides of $K_6$ are joined by
\begin{multline*}
S(S(a)S(b)S(S(c)d)e)f
\sim_S S(S(S(a)S(b)S(c))de)f
\sim_S S(S(a)S(b)S(S(c)de))f\\
\sim_S S(a)S(b)S(S(S(c)de)f)
\sim_S S(a)S(S(b)S(S(c)de))f
\sim_S S(a)S(S(S(b)S(c))de)f\\
\sim_S S(a)S(S(b)S(S(c)d)e)f
\sim_S S(S(a)S(b))S(S(c)d)ef\\
\sim_S S(a)S(S(b)S(S(c)d))ef
\sim_S S(a)S(b)S(S(S(c)d)e)f.
\end{multline*}
Similarly, the two sides of $E_6$ are joined by
\begin{multline*}
S(S(a)S(S(b)c)S(d)e)f
\sim_S S(S(S(a)S(b))cS(d)e)f
\sim_S S(S(a)S(S(b)cS(d))e)f\\
\sim_S S(S(a)S(b)cS(S(d)e))f
\sim_S S(a)S(b)cS(S(S(d)e)f)
\sim_S S(a)S(S(b)cS(S(d)e))f\\
\sim_S S(S(a)S(b))cS(S(d)e)f
\sim_S S(a)S(S(b)c)S(S(d)e)f\\
\sim_S S(a)S(S(S(b)c)S(d))ef
\sim_S S(S(a)S(S(b)c))S(d)ef
\sim_S S(a)S(S(S(b)c)S(d)e)f.
\end{multline*}
Together with Theorem~\ref{thm:validity}, these paths show that, for
\[
 I^{\succ}=\Id_{\succ}\langle F_r,G_r\ (r\ge3),K_6,E_6\rangle,
\]
we have $I^{\succ}\subseteq\Ker\psi$.

For a pure $\succ$-monomial, put $D_j=\sum_{i\leq j}d_i$ and
$H_j=j-D_j$, with $H_0=0$. Lemmas~\ref{lem:pure-vectors} and
\ref{lem:center-tree} give $H_j>0$ for $j<n$ and
$M(T)=\{i<n:d_i(T)=0\}$.

\begin{lemma}
\label{pure:lem:prefix}
For a pure $\succ$-monomial $T$ and $t<n$,
\begin{equation*}
[1,t]\in\Ical(\psi(T))\quad\Longleftrightarrow\quad
D_t=t-1\quad\Longleftrightarrow\quad H_t=1.
\end{equation*}
\end{lemma}
\begin{Proof}
In the reconstruction of Lemma~\ref{lem:pure-vectors}, read from left
to right for $\succ$, $H_t$ is the number of monomials remaining after
the first $t$ leaves. If $H_t=1$, these leaves form a submonomial $P$.
Since $t<n$, it is later used as a left factor, so $S(\psi(P))$ occurs.
Conversely, an interval $[1,t]$ comes from a left submonomial on these
leaves, so they form one component and $H_t=1$.
\end{Proof}

For $s\geq2$, write
$R_s(A_1,\ldots,A_s;X)=A_1\succ(A_2\succ\cdots\succ(A_s\succ X))$.

\begin{lemma}
\label{pure:lem:Ts}
If $X=(C\succ D)\succ E$, then for every $s\geq2$,
\begin{equation}
R_s(A_1,\ldots,A_s;X)\succ Y
\equiv R_s(A_1,\ldots,A_s;X\succ Y)\pmod{I^{\succ}}.
\label{pure:Ts}
\end{equation}
\end{lemma}
\begin{Proof}
For $s=2$, this is $K_6$. Suppose it holds for $s$. In the left side
for $s+1$, use $F_3$ to replace
$A_s\succ(A_{s+1}\succ((C\succ D)\succ E))$ by
$A_s\succ((A_{s+1}\succ(C\succ D))\succ E)$.
Apply the induction hypothesis to the first $s$ entries with
$X'=(A_{s+1}\succ(C\succ D))\succ E$, and then apply $F_4$ backwards
to the inner expression. The result is
$R_{s+1}(A_1,\ldots,A_{s+1};X\succ Y)$.
\end{Proof}

\begin{remark}
\label{pure:rem:degree-effects}
Every pure $\succ$-submonomial has center at its last leaf. Thus a
replacement on the same support leaves all ancestor centers unchanged.
The change in the degree vector can be read directly from the identity.
Applied backwards, $F_r$ moves one unit from the last leaf of $c_2$ to
that of $c_r$, and $G_r$ moves one unit from the last leaf of $b_3$ to
that of $d$. Applied forwards, $E_6$ moves one unit from the last leaf
of $e$ to that of $f$, while \eqref{pure:Ts} moves $s$ units from the
last leaf of $E$ to that of $Y$. Each move preserves the zero set and
moves the indicated units to the right.
\end{remark}

For a singleton set $M\subseteq\{1,\ldots,n-1\}$ containing $1$, let
\begin{equation*}
r=\max\{j:\{1,\ldots,j\}\subseteq M\}.
\end{equation*}

\begin{lemma}
\label{pure:lem:r1-prefix}
Suppose $r=1$ and $q=\min(M\setminus\{1\})$ exists. Then
$d_2=\cdots=d_{q-1}=1$, and the first $q-1$ leaves form the literal
left comb $L_{q-1}(x_1,\ldots,x_{q-1})$.
\end{lemma}
\begin{Proof}
For $2\leq i<q$, we have $d_i>0$ and $D_i<i$. If
$d_2=\cdots=d_{i-1}=1$, then $D_i=i-2+d_i<i$, so $d_i=1$.
Induction proves the degree assertion. The reconstruction gives the
stated left comb.
\end{Proof}

We contract such a prefix before applying the next lemma. Thus it
suffices there to assume $r\geq2$. Call $d_t>1$ an excess coordinate.

\begin{lemma}
\label{pure:lem:leftmost-excess}
Suppose $r\geq2$, and let $t<n$ be the first excess coordinate.
Put $h=H_{t-1}$ and $k=d_t$. Then $2\leq k\leq h$.
Except for the locked configurations described in the proof, an
instance of $F_r,G_r,E_6$ or \eqref{pure:Ts} moves at least one excess
unit from position $t$ to a later position.
\end{lemma}
\begin{Proof}
For $i<t$, all degrees are zero or one. Immediately before processing
$x_t$, the reconstruction therefore has components $C_1,\ldots,C_h$,
each a left comb on consecutive leaves. A zero starts a new component;
a one replaces the last component $C$ by $C\succ x_i$. Since $r\geq2$
and no excess has occurred, $C_1=x_1$.

Suppose first that $k<h$. After the $k$ attachments at $x_t$, the last
component is
\begin{equation*}
X=C_{h-k+1}\succ(C_{h-k+2}\succ\cdots\succ(C_h\succ x_t)),
\end{equation*}
and the component immediately before it is $a=C_{h-k}$.
Consider the first subsequent attachment joining $a$ to a component
containing $X$. Before this happens, the latter has the form
$L_{m+1}(X,V_1,\ldots,V_m)$ with $m\geq1$: at least one nonempty
block must be processed before two existing components can be joined.
Since $k\geq2$, write $X=b\succ(c_1\succ c_2)$, where
$b=C_{h-k+1}$, $c_1=C_{h-k+2}$ and $c_2$ ends at $x_t$.
The attachment contains
\[
a\succ L_{m+1}\bigl(b\succ(c_1\succ c_2),V_1,\ldots,V_m\bigr).
\]
Applying $F_{m+2}$ backwards moves one unit from position $t$ to the
last leaf of $V_m$.

Now suppose $k=h$. All components are joined into
\begin{equation*}
X=C_1\succ(C_2\succ\cdots\succ(C_h\succ x_t)).
\end{equation*}
If every $C_i$ is a leaf, then $t=r+1$ and $d_{r+1}=r$.
We retain this configuration and call it a locked first return.
Otherwise, let $p$ be the last index for which $C_p$ is nontrivial.
Since $C_1=x_1$, we have $p\geq2$.

If $p\geq3$, write $C_p=U\succ V$ and
$E_0=C_{p+1}\succ(\cdots\succ(C_h\succ x_t))$, with $E_0=x_t$
when $p=h$. Then
\[
X=R_{p-1}(C_1,\ldots,C_{p-1};(U\succ V)\succ E_0).
\]
Since $t<n$, this component is later a left factor of $X\succ Y$.
Relation \eqref{pure:Ts} moves $p-1$ units from position $t$ to the
last leaf of $Y$.

If $p=2$ and $h\geq3$, write $C_2=U\succ V$. The remaining right
comb is nontrivial, say $D\succ E$, so
$X=C_1\succ((U\succ V)\succ(D\succ E))$.
At its first attachment to a right factor $Y$, apply $E_6$ to
$X\succ Y$. This moves one unit from $t$ to the last leaf of $Y$.

It remains that $p=h=2$. Before $t$, the height equals the number
of zero coordinates, so $M\cap\{1,\ldots,t-1\}=\{1,2\}$.
Moreover, $X=C_1\succ(C_2\succ x_t)$ with $C_2$ a nontrivial left
comb. It is a substitution into $T_4$, with the fourth argument ending
at $t$. Subsequent coordinates equal to one extend it by
$T_j\mapsto T_{j+1}$.
If a later singleton occurs, consider the first attachment joining its
component to this one. The right component is nontrivial, say
$C\succ D$, and the attachment is $T_j\succ(C\succ D)$.
Applying $G_{j-1}$ backwards moves the excess from $t$ to the right.
If no later singleton occurs, then $M=\{1,2\}$, the terminal pair
sector of Corollary~\ref{cor:pair-zero}.
\end{Proof}

Put
\begin{equation*}
\Omega(T)=\sum_{i=1}^{n-1}(n-i)(d_i(T)-1)_+.
\end{equation*}

\begin{lemma}
\label{pure:lem:termination}
Every move in Lemma~\ref{pure:lem:leftmost-excess} strictly decreases
$\Omega$.
\end{lemma}
\begin{Proof}
Each moved excess unit has weight $n-i$ before the move and at most
$n-j$ afterwards, where $j>i$. No excess is created to its left.
Hence $\Omega$ decreases. Since it is a nonnegative integer, the
process terminates.
\end{Proof}

\noindent\textbf{Canonical pure monomials.}

We define $\Ccal_n(M)$ recursively, starting with
$\Ccal_1(\varnothing)=\{x_1\}$. For $n\geq2$ and $M=\{1\}$, the only
vector is $(0,1,\ldots,1)$. For $M=\{1,2\}$, the vectors are
\begin{equation*}
d_i=
\begin{cases}
0,&i=1,2,\\
2,&i=t+1,\\
1,&\text{otherwise},
\end{cases}
\qquad 2\leq t\leq n-1.
\end{equation*}
The corresponding monomials are $P_n^t$.

Suppose $M=\{1,\ldots,r\}$ with $r\geq3$. We take the locked form
\begin{equation*}
L_M:\quad d=(\underbrace{0,\ldots,0}_{r},r,1,\ldots,1)
\end{equation*}
and the mobile form $B_M$, defined for any $M$ by
\begin{equation*}
d_i(B_M)=
\begin{cases}
0,&i\in M,\\
1,&i<n,\ i\notin M,\\
|M|,&i=n.
\end{cases}
\end{equation*}
These coincide when $n=r+1$; otherwise they give two forms.

Suppose now that a later singleton
$q=\min(M\setminus\{1,\ldots,r\})$ exists. If $r=1$, let $P$ be
the forced left comb on the first $q-1$ leaves. If $r\geq2$, let $P$
be the prefix with $d_{r+1}=r$ and
$d_{r+2}=\cdots=d_{q-1}=1$. Contracting $P$ gives arity
$n'=n-q+2$ and singleton set
\begin{equation*}
M'=\{1,2\}\cup\{m-q+2:m\in M,\ m>q\}.
\end{equation*}
Write $\operatorname{Exp}_P(\overline T)$ for substitution of $P$ at
the first leaf, with the original leaf labels restored. Set
\[
\Ccal_n(M)=
\begin{cases}
\operatorname{Exp}_P(\Ccal_{n'}(M')),&r=1,\\
\{B_M\}\sqcup\operatorname{Exp}_P(\Ccal_{n'}(M')),&r\geq2.
\end{cases}
\]
Together with the preceding terminal cases, this defines the canonical
set, since $n'<n$.

\begin{proposition}
\label{pure:prop:normalization}
Every pure $\succ$-monomial with singleton set $M$ is congruent modulo
$I^{\succ}$ to an element of $\Ccal_n(M)$.
\end{proposition}
\begin{Proof}
Proceed by induction on $n$. The one- and two-singleton terminal cases
already have the stated forms. If $r=1$ and a later singleton exists,
Lemma~\ref{pure:lem:r1-prefix} gives the literal prefix $P$.
Contract it, apply induction in arity $n'$, and substitute $P$ back.
This preserves congruence because $I^{\succ}$ is an operadic ideal.

Suppose $r\geq2$. Apply Lemma~\ref{pure:lem:leftmost-excess}
repeatedly to the first excess coordinate. By
Lemma~\ref{pure:lem:termination}, either all excess reaches the last
coordinate, giving $B_M$, or a locked first return occurs. In the latter
case, contract the prefix through $q-1$ if a later singleton $q$ exists,
apply induction, and expand. If no later singleton exists, the locked
form is already one of the terminal forms. These are exactly the
possibilities in the definition of $\Ccal_n(M)$.
\end{Proof}

\begin{lemma}
\label{pure:lem:literal-rigid}
Every forced or locked prefix used above has an operated image to which
\eqref{R} cannot be applied internally in either direction.
\end{lemma}
\begin{Proof}
For a forced prefix, the operated images are
$p_2=S(x_1)x_2$ and $p_{j+1}=S(p_j)x_{j+1}$.
At each new level, $S(p_j)$ is first, so the forward relation has no
nonempty factor before it. Since $p_j$ ends in a variable, the reverse
relation cannot apply either. Induction proves the assertion.

A locked prefix starts with the right comb on $r+1$ leaves, whose image
is $S(x_1)\cdots S(x_r)x_{r+1}$. Each $S$-argument is a variable.
Further left extensions again have the form $p_{j+1}=S(p_j)x_{j+1}$,
so the same induction applies.
\end{Proof}

\begin{lemma}
\label{pure:lem:prefix-persistent}
The interval $[1,q-1]$ of a forced or locked prefix remains present
along every $S$-path from its canonical representative.
\end{lemma}
\begin{Proof}
In the forced case, start with the invariant interval $[1,1]$.
In the locked case, the initial rigid block gives $[1,r+1]$.
Its first removal would require a pivot $[p,r+1]$ with $p>1$, which
is absent from this block.

Suppose $[1,s-1]$ persists and $s<q$ is a nonsingleton position.
Changing $[1,s]$ requires a pivot $[p,s]$ with $p>1$.
For $p<s$, this crosses $[1,s-1]$. For $p=s$, it is the singleton
$[s,s]$, which is absent. Thus $[1,s]$ also persists, and induction
reaches $q-1$.
\end{Proof}

\begin{lemma}
\label{pure:lem:S-contraction}
Let $P$ be a forced or locked prefix. If
$\psi(\operatorname{Exp}_P(U))\sim_S
\psi(\operatorname{Exp}_P(V))$, then $\psi(U)\sim_S\psi(V)$.
\end{lemma}
\begin{Proof}
By Lemmas~\ref{pure:lem:literal-rigid} and
\ref{pure:lem:prefix-persistent}, the argument $\psi(P)$ is internally
rigid and its enclosing interval persists. Every move is therefore
disjoint from this block or uses it as a whole. Replacing $\psi(P)$
by one variable throughout the path gives the required path between
$\psi(U)$ and $\psi(V)$.
\end{Proof}

\begin{lemma}
\label{pure:lem:mobile-locked}
If the initial singleton run has length $r\geq2$ and a later singleton
$q$ exists, then $B_M$ is not $S$-equivalent to any locked canonical
expansion.
\end{lemma}
\begin{Proof}
For $B_M$, we have $D_{q-1}=q-1-r<q-2$. By
Lemma~\ref{pure:lem:prefix}, its image does not contain $[1,q-1]$.
Every locked expansion contains this interval, which persists by
Lemma~\ref{pure:lem:prefix-persistent}. The classes are therefore distinct.
\end{Proof}

\begin{proposition}
\label{pure:prop:canonical-separation}
Distinct elements of $\Ccal_n(M)$ have images in distinct $S$-classes.
\end{proposition}
\begin{Proof}
Proceed by induction on $n$. There is one form for $M=\{1\}$, and
Corollary~\ref{cor:pair-zero} separates the forms for $M=\{1,2\}$.
When a later singleton exists, Lemma~\ref{pure:lem:mobile-locked}
separates the mobile form from the locked expansions. Two distinct
expansions contract to distinct canonical forms in smaller arity, so
Lemma~\ref{pure:lem:S-contraction} and induction separate them.
For $r=1$, only this contraction argument is needed.

If there is no later singleton and $r\geq3$, the two forms coincide
when $n=r+1$. Otherwise the locked form has the persistent interval
$[1,r+1]$, while $D_{r+1}(B_M)=1<r$. Thus
Lemma~\ref{pure:lem:prefix} separates them as well.
\end{Proof}

Let $\Mag(n,M)$ be the span of pure $\succ$-monomials of arity $n$
with singleton set $M$. Write
$I^{\succ}(n,M)=I^{\succ}(n)\cap\Mag(n,M)$.

\begin{lemma}
\label{pure:lem:sector-ideal}
For every $n$, $I^{\succ}(n)=\bigoplus_M I^{\succ}(n,M)$.
\end{lemma}
\begin{Proof}
Every generating binomial of $I^{\succ}$ has $S$-equivalent images,
and this remains true under substitution. By
Lemma~\ref{lem:singleton}, each such consequence belongs to one
singleton sector. These consequences span the ideal.
\end{Proof}

\begin{theorem}
\label{pure:thm:pure-complete}
The ideal of identities of $a\succ b=S(a)b$ is
\begin{equation*}
\ker\psi=\Id_{\succ}\langle F_r,G_r\ (r\geq3),K_6,E_6\rangle.
\end{equation*}
\end{theorem}
\begin{Proof}
We have already proved $I^{\succ}\subseteq\ker\psi$.
For a fixed $n,M$, Proposition~\ref{pure:prop:normalization} shows
that $\Ccal_n(M)$ spans $\Mag(n,M)/I^{\succ}(n,M)$.
By Proposition~\ref{pure:prop:canonical-separation}, its images lie
in distinct operated-word classes and hence are linearly independent.
The induced surjection from this quotient to $\operatorname{Im}\psi(n,M)$
is therefore an isomorphism. Thus
$I^{\succ}(n,M)=\ker\psi(n,M)$. Summing over $M$ proves the assertion.
\end{Proof}

\section{The upper bound}
\label{sec:small}

Write $P\equiv_n Q$ if $P-Q\in D_n$. In particular, any instance
of a valid identity of arity less than $n$ has this property, as does
an identity applied in a proper binary submonomial.

We shall use the identities
\begin{align}
(a\prec b)\succ(c\succ d)&=((a\prec b)\prec c)\succ d,\label{chi}\\
(a\succ(b\succ c))\prec(d\succ e)
&=a\succ((b\succ c)\succ(d\succ e)).\label{Cfive}
\end{align}
For \eqref{chi}, the operated images satisfy
\[
S(aS(b))S(c)d\sim_S aS(S(b)S(c))d
\sim_S aS(b)S(S(c)d)\sim_S S(aS(b)S(c))d.
\]
For \eqref{Cfive}, read the following path in reverse:
$$S(a)S(S(b)c)S(d)e\sim_S S(S(a)S(b))cS(d)e\sim_S S(a)S(S(b)cS(d))e\sim_S S(a)S(b)cS(S(d)e).$$
Thus every arity-$n$ instance of these identities belongs to $D_n$
when $n\geq8$.

\begin{lemma}
\label{lem:switch}
Modulo \eqref{NB} and \eqref{chi}, $U\prec V=U\succ V$ whenever
$U$ contains $\prec$ and $V$ contains $\succ$.
Modulo \eqref{NB} and \eqref{Cfive},
\begin{equation}
a\succ(B\prec V)=a\succ(B\succ V)
\label{Cgate}
\end{equation}
whenever both $B$ and $V$ contain $\succ$.
\end{lemma}
\begin{Proof}
By Lemma~\ref{lem:rootfactor}, write $U=u\prec v$ and $V=w\succ z$.
Identity \eqref{chi}, followed by
$(U\prec w)\succ z=U\prec(w\succ z)$, proves the first assertion.
For the second, write $B=b\succ c$ and $V=d\succ e$, and apply
\eqref{Cfive} to $a\succ(B\prec V)=(a\succ B)\prec V$.
\end{Proof}

Suppose $V$ contains $\succ$. Then $\RR(\phi(V))$ holds, and
\eqref{generalB} and \eqref{generalE} give
\begin{align}
U\prec(z\prec V)&=(U\succ z)\prec V,
\qquad \LL(\phi(U)),\label{Bold}\\
(a\succ(b\succ c))\prec(d\prec V)
&=(a\succ((b\succ c)\succ d))\prec V.\label{Eold}
\end{align}
We also need that these relations are consequences of smaller ones.
Write $V=C\succ D$. The two sides of \eqref{Bold} become
$(U\prec(z\prec C))\succ D$ and $U\succ((z\prec C)\succ D)$.
They are obtained from the valid relation
$(U\prec w)\succ t=U\succ(w\succ t)$, which follows from
\eqref{H0}. This relation has arity $\arity(U)+2$, strictly smaller
than the given instance, since $w=z\prec C$ has at least two leaves.
Similarly, \eqref{Eold} is an instance of a valid six-input identity
when $V=C\succ D$, so its arity-$n$ instances belong to $D_n$ for
$n\geq8$.

\begin{lemma}
\label{lem:flat}
If pure $\succ$-monomials $P,Q$ have the same ordered leaves and
singleton set, then $A\prec P=A\prec Q$ modulo \eqref{NB}, \eqref{chi}.
The reflected assertion holds for pure $\prec$-monomials in $P\succ A$.
\end{lemma}
\begin{Proof}
Apply $A\prec(B\succ C)=(A\prec B)\succ C$ down the left spine
until the first factor is $A\prec x_1$. It contains $\prec$.
For every remaining nontrivial pure right factor, switch the root by
Lemma~\ref{lem:switch} and repeat. At each step the arity of the
unprocessed right factors decreases.

The result is a left comb starting at $A\prec x_1$. The leaf $x_i$
is adjoined by $\prec$ when $i\in M(P)$ and by $\succ$ otherwise.
Indeed, in a right factor its left child is adjoined by $\prec$ and
its right child by $\succ$; apply the same rule recursively to a
nontrivial left child. Thus the resulting comb depends only on the
singleton set, proving the assertion.
\end{Proof}

\noindent\textbf{Interior sectors.}
\label{sec:interior}

\begin{lemma}
\label{lem:expose}
If $1,n\notin M(T)$, then modulo \eqref{NB}, \eqref{chi}, the monomial
$T$ can be written as $x_1\prec V$ and as $U\succ x_n$, where $V$
contains $\succ$ and $U$ contains $\prec$.
More generally, every mixed monomial with $n\notin M(T)$ can be written
as $P\prec V$, where $P$ is pure $\succ$ and $V$ contains $\succ$.
\end{lemma}
\begin{Proof}
For the first assertion, move a $\succ$ to the root and write
$T=U\succ V$. Since the first leaf is not a singleton, $U$ is
nontrivial and contains $\prec$. Write $U=A\prec B$ to obtain
$A\prec(B\succ V)$. If $A$ still contains $\prec$, switch this root
to $\succ$, factor $A=A_1\prec A_2$, and obtain
$A_1\prec(A_2\succ(B\succ V))$.
The left factor becomes shorter each time. It must end at $x_1$,
since a nontrivial pure $\succ$ prefix would make the first leaf a
singleton. Reflection gives $U\succ x_n$.

For the general assertion, move $\prec$ to the root. The last leaf
condition makes the right factor nontrivial and forces it to contain
$\succ$. The same procedure stops when the left factor is pure $\succ$.
\end{Proof}

\begin{theorem}
\label{thm:interior}
If $1,n\notin M$, all normal monomials in $\Ncal_n(M)$ are equal
modulo \eqref{NB}, \eqref{chi}. Hence $K_n(M)\subseteq D_n$ for $n>4$.
\end{theorem}
\begin{Proof}
Proceed by induction on $n$. The first case, $n=3$ and $M=\{2\}$,
has one normal monomial. If $n-1\notin M$, write $T=U\succ x_n$;
then $U$ is again interior. If $2\notin M$, write $T=x_1\prec V$;
then $V$ is interior. Apply induction in the appropriate submonomial.

It remains that $2,n-1\in M$ and $n\geq4$. Write $T=U\succ x_n$
and normalize $U$. Its last leaf is a singleton, so $U=A\prec P$
with $P$ pure $\prec$. The factor $A$ contains $\prec$. Otherwise,
the first leaf condition forces $A=x_1$, and a nontrivial $P$ would
make leaf $2$ a nonsingleton. If $P$ were a leaf, we would have $n=3$.
Thus $T=A\prec(P\succ x_n)$.

By the reflected form of Lemma~\ref{lem:flat}, write $P\succ x_n$
as a right comb ending in $Y=x_{n-1}\succ x_n$. Move $A$ through
this comb. At $C\,\mu\,(z\,\nu\,W)$, the left factor contains
$\prec$ and the right factor contains $\succ$, since it contains $Y$.
If necessary, switch $\mu$ to the operation different from $\nu$,
and apply \eqref{NB} to obtain $(C\,\mu\,z)\,\nu\,W$.
The remaining right comb becomes shorter. At its end, switch the root
to obtain $C'\succ Y$.

Contract $Y$ to one leaf. The resulting monomial is interior, has
arity $n-1$, and has singleton set $M\setminus\{n-1\}$, with labels
shifted after contraction. Induction gives a common representative.
Substituting $Y$ back proves the assertion, since the same choice of
case applies to every monomial of the fixed sector.
\end{Proof}

\noindent\textbf{Rigid ends.}
\label{sec:pure-use}

\begin{lemma}
\label{lem:old-rigid}
Suppose $\phi(T)$ contains a rigid block $S(u)$ on a proper global
prefix of at least two leaves, where $u$ is regular and ends in a
direct variable. Then every fibre binomial $T-T'$ belongs to $D_n$.
The reflected assertion holds for a rigid suffix whose argument starts
in a direct variable.
\end{lemma}
\begin{Proof}
By Lemma~\ref{lem:rigid-boundary}, the block is unchanged along a path
to $\phi(T')$. Contract its argument to one variable. The endpoints
remain regular, and their preimages give a kernel relation of arity
$n-\arity(u)+1<n$. Expanding by the same preimage of $u$ recovers
$T,T'$ modulo \eqref{NB}, by Lemma~\ref{lem:literal-unique}.
Their difference is therefore a consequence of this smaller relation.
\end{Proof}

It remains to consider fibres without nontrivial rigid proper ends.
This condition is preserved along an $S$-path, since persistence also
applies to the reversed path.

\begin{lemma}
\label{lem:pure-no-prefix}
Let $P$ be pure $\succ$ with no nontrivial rigid proper prefix.
If its initial singleton run has length one, then $P$ has at most two
leaves. If the run has length at least two, its canonical form is
$B_M$, or $P_m^t$ in the pair sector, with $t\geq3$ unless $m=3$.
If $\LL(\phi(P))$ fails, then $P$ is a literal right comb.
\end{lemma}
\begin{Proof}
Every forced or locked canonical form excluded by the statement has
a persistent rigid proper prefix. The left comb with singleton set
$\{1\}$ also has one when its arity is at least three, as does $P_m^2$
when $m>3$. The canonical list therefore gives the first assertions.

For the last, write $P=A\succ B$. By Lemma~\ref{lem:hiddenarm},
$\LL(\phi(A))$ fails and $B$ is a right comb. If $A$ were nontrivial,
its internally irreducible image would give the forbidden prefix
$S(\phi(A))$. Thus $A$ is a leaf and $P$ is a right comb.
\end{Proof}

\begin{lemma}
\label{lem:restricted-cancel}
For pure $\succ$-monomials $P,Q$,
\[
\phi(P\prec z)\sim_S\phi(Q\prec z)
\quad\Longrightarrow\quad \phi(P)\sim_S\phi(Q).
\]
The reflected assertion holds for $z\succ P$ with $P$ pure $\prec$.
\end{lemma}
\begin{Proof}
Replace $P,Q$ by pure canonical forms. Their singleton sets agree.
Proceed by induction on their arity. A mobile form has no nontrivial
rigid proper prefix, whereas a forced or locked expansion does.
Adjoining $\prec z$ preserves these prefixes. By persistence, it
cannot make a mobile form equivalent to a different forced or locked
form. Two forced or locked forms have the same prescribed prefix;
contract it along the path and apply induction.

There is one form in the one-singleton terminal sector. In the pair
sector, the monomials $P_m^t\prec z$ belong to the zero-kernel sector
$\{1,2,m+1\}$, so different $t$ remain distinct. If the initial run
has length $r\geq3$ and there is no later singleton, the two terminal
forms are separated by the rigid prefix on $r+1$ leaves, unless they
coincide. These are all canonical cases. Reflection proves the other
assertion.
\end{Proof}

\noindent\textbf{One boundary singleton.}
\label{sec:one-boundary}

\begin{proposition}
\label{prop:one-boundary}
Suppose $1\in M$, $n\notin M$, and there is no nontrivial rigid proper
prefix. If the initial singleton run has length at least two, every
monomial in the sector is equivalent to a pure $\succ$-monomial.
For $n\geq8$, this equivalence holds modulo $D_n$.
If the initial run has length one, every monomial is equal modulo
\eqref{NB}, \eqref{chi} to $x_1\succ W$, where $W$ is interior or,
in arity two, a leaf. All monomials with this singleton set are then
equal modulo these relations.
\end{proposition}
\begin{Proof}
For a mixed monomial, Lemma~\ref{lem:expose} gives $P\prec V$ with
$P$ pure $\succ$ and $V$ containing $\succ$. The first leaf condition
makes $P$ nontrivial. The absence of a rigid prefix is preserved by
these transformations.

We may replace $V$ by any monomial with its singleton set and the same
nonsingleton last leaf. Indeed, Theorem~\ref{thm:interior} identifies
the corresponding monomials $y\prec V$ for a new variable $y$, and we
then substitute $P$ for $y$. Repeated boundary exposure gives
$V=z_1\prec(z_2\prec\cdots\prec V_0)$, where the $z_i$ are leaves
and the first leaf of $V_0$ is a singleton. This process stops because
every nontrivial regular word has a singleton. We may choose $V_0$
pure $\succ$ with the prescribed singleton set; the vector of $B_M$
gives such a monomial.

Suppose $\LL(\phi(P))$ holds. Apply \eqref{Bold} successively to
adjoin the $z_i$ to $P$ by left $\succ$ extensions. The property
$\LL$ is preserved. By Lemma~\ref{lem:hiddenarm}, reduce the resulting
prefix modulo $F_3$ to $T_k(U_1,\ldots,U_k)$, with all $U_i$ pure.
Apply \eqref{Bold} backwards to its final extensions until the left
factor is $T_4$. The result is
\[
T_4(U_1,U_2,U_3,U_4)\prec
\bigl(U_5\prec(\cdots\prec V_0)\bigr).
\]
By \eqref{Cgate}, this equals
\[
U_1\succ\bigl(((U_2\succ U_3)\succ U_4)\succ
(U_5\prec(\cdots\prec V_0))\bigr).
\]
At every remaining $\prec$ in the tail, first apply
$B\succ(z\prec W)=(B\succ z)\prec W$ and then \eqref{Cgate} inside
the outer $U_1\succ(\cdot)$. Both required factors contain $\succ$:
one is $B\succ z$, and the other contains $V_0$.
Each step removes one $\prec$ from the tail. We obtain the pure monomial
$U_1\succ(L_{k-1}(U_2,\ldots,U_k)\succ V_0)$.

Suppose $\LL(\phi(P))$ fails. By Lemma~\ref{lem:pure-no-prefix},
$P$ is a right comb. If it has at least three leaves, write
$P=a\succ(b\succ c)$. When a first $z_1$ occurs, \eqref{Eold}
replaces the prefix by $a\succ((b\succ c)\succ z_1)$, for which
$\LL$ holds, and the preceding argument applies. If there are no $z_i$,
\eqref{Cfive} gives a pure monomial directly.

The remaining prefix is $P=x_1\succ x_2$. Then
$P\prec V=x_1\succ(x_2\prec V)$, with the inner monomial interior,
and the initial singleton run has length one. Conversely,
Lemma~\ref{lem:pure-no-prefix} forces this prefix when the run has
length one. Theorem~\ref{thm:interior} identifies all the inner
monomials with the prescribed singleton set.

For $n\geq8$, every use of \eqref{chi}, \eqref{Cfive} and \eqref{Eold}
belongs to $D_n$. The same holds for \eqref{Bold} by the arity
calculation above. Pure changes in the proper prefix also belong to
$D_n$. This proves the stated congruence.
\end{Proof}

Thus every one-boundary kernel class reduces modulo $D_n$ to a pure
kernel class. By Theorem~\ref{pure:thm:pure-complete}, for $n\geq8$
these classes are spanned by $F_{n-2},G_{n-3}$. Every smaller pure
generator, including $K_6,E_6$, contributes only to $D_n$.
Reflection gives $\rho F_{n-2},\rho G_{n-3}$ at the other boundary.

\noindent\textbf{Two boundary singletons.}
\label{sec:both}

\begin{lemma}
\label{lem:short-run}
If $1,n\in M$ and either boundary singleton run has length one, then
$K_n(M)\subseteq D_n$ for $n\geq8$.
\end{lemma}
\begin{Proof}
By Lemma~\ref{lem:old-rigid}, assume there are no rigid proper ends.
Apply reflection if necessary so that the initial run has length one.
A normal monomial ending at a singleton has a final left spine of
$\prec$-vertices, so write
\[
T=((U\prec Q_1)\prec\cdots)\prec Q_j,
\]
where each $Q_i$ is pure $\prec$ and the last leaf of $U$ is a
nonsingleton. Its first leaf is a singleton, so $U$ is not a leaf.
It has no rigid proper prefix. Proposition~\ref{prop:one-boundary}
gives $U=x_1\succ W$, with $W$ interior or a leaf. Move the final
$\prec$ attachments inside by \eqref{NB} to obtain $x_1\succ X$.
The first leaf of $X$ is a nonsingleton, its last is a singleton, and
it has no rigid suffix.

If the terminal run has length at least two, apply the reflected
one-boundary reduction inside $X$ to obtain $x_1\succ Q$, with $Q$
pure $\prec$. For two endpoints in the same fibre,
Lemma~\ref{lem:restricted-cancel} gives $\phi(Q)\sim_S\phi(Q')$.
Their difference is therefore a consequence of a proper pure identity.
If the terminal run has length one, the reduction gives
$x_1\succ(Y\prec x_n)$, with $Y$ interior or a leaf, and
Theorem~\ref{thm:interior} identifies these monomials.
All changes inside $U$ and $X$ occur in proper submonomials, so they
belong to $D_n$ in either case.
\end{Proof}

For pure $\succ$ $P$ and pure $\prec$ $Q$, let
$P\star Q=P[\text{last leaf}:=Q]$. The first relation in \eqref{NB}
gives
\begin{equation}
P\star Q=Q[\text{first leaf}:=P].
\label{graft}
\end{equation}
Indeed, it moves a $\prec$ attachment down the right spine of $P$,
one vertex at a time. The right side is normal, since $P$ is inserted
on the left spine of $Q$. This relation preserves centers and degree
vectors, including in larger contexts.

\begin{proposition}
\label{prop:concentrated}
Suppose $M\supseteq\{1,2,n-1,n\}$ and there are no nontrivial rigid
proper ends. Every normal monomial is equivalent to a normal monomial
$A(M,k)$, $k\notin M$, with
\begin{equation}
d_i=0\ (i\in M),\qquad
 d_i=1\ (i\notin M,\ i\ne k),\qquad d_k=|M|.
\label{concentrated-vector}
\end{equation}
For $n\geq8$, the difference belongs to $D_n$ plus the span of the
arity-$n$ corner identities $J,H,\rho J$. If $|M|>4$, it belongs to
$D_n$.
\end{proposition}
\begin{Proof}
Use the final $\prec$ spine from the preceding proof. Its bottom factor
has initial singleton run of length at least two. By
Proposition~\ref{prop:one-boundary}, applied in this proper factor,
it may be made pure $\succ$. Thus the monomial is a graft $P\star Q$
with common leaf $k$, where $3\leq k\leq n-2$.
Normalize both proper factors by the pure theorem. A rigid proper
prefix of $P$ or suffix of $Q$ would remain a rigid end of the whole
graft. Hence each factor is mobile or belongs to its terminal pair
list $P_m^t$ or its reflection.

Consider a left pair factor $P_k^t$. For $k=3$ it is already mobile;
otherwise $t\geq3$. If the right factor has property $\RR$, reduce
it internally by the reflected $F_3$ to
$Q_\ell(V_1,\ldots,V_\ell)$ with $\ell\geq4$.
If it does not, the reflected form of Lemma~\ref{lem:pure-no-prefix}
makes it a literal left $\prec$-comb. Its length is at least three,
so write it as $Q_3(V_1,V_2,V_3)$, where
$Q_3=(x_1\prec x_2)\prec x_3$.

Compare $P_k^t\star Q_\ell$ and $P_k^{k-1}\star Q_\ell$ before
substitution. If $t=k-1$, no change is needed. Otherwise $k\geq5$,
and these monomials have arity $N=k+\ell-1\geq7$. Put $K=N-4$.
Their singleton set is the four-singleton corner. For the first
monomial, the three excess positions are
\[
\begin{cases}
(t+1,k,k),&\ell=3,\\
(t+1,k,N-3),&\ell\geq4.
\end{cases}
\]
The left factor contributes its excess at $t+1$, the right factor
contributes its own excess, and their two positive degrees at $k$
contribute one further excess there. By \eqref{corner-excess}, the
corner coordinates are
\[
\begin{cases}
(k-t,K,K+1),&\ell=3,\\
(k-t,k-2,K),&\ell\geq4.
\end{cases}
\]
For $3\leq t\leq k-1$, these belong to $\mathcal M_K$:
in the first case $1\leq k-t\leq K-1$, and in the second
$1\leq k-t<k-2<K$. Proposition~\ref{prop:corner-mobile} therefore
identifies the two monomials. Substitute the $V_i$, including $V_1$
at the common leaf, and undo the internal changes in the right factor.
This makes the left pair factor mobile. Reflection and \eqref{graft}
give the same construction for a right pair factor.

If a substituted block is nontrivial, the relation used has arity
strictly less than $n$, so its consequence belongs to $D_n$.
If it has arity $n$, all substituted blocks are leaves and $|M|=4$;
Corollary~\ref{cor:corner-upper} expresses the difference by the stated
corner identities modulo $D_n$.

Both factors are now mobile. In the normal expression on the right
of \eqref{graft}, their degrees add at $k$ and are unchanged elsewhere.
Their singleton sets partition $M$, giving \eqref{concentrated-vector}.
These mobile grafts also show that $A(M,k)$ has no rigid proper ends.
\end{Proof}

\noindent\textbf{Concentrated degree vectors.}
\label{sec:concentrated-induction}

\begin{lemma}
\label{lem:extreme-cherry}
Suppose the first positive coordinate of a normal monomial is at
$\alpha>1$. Then it contains the literal submonomial
$x_{\alpha-1}\succ x_\alpha$. Contracting it replaces
$(0,d_\alpha)$ by $d_\alpha-1$ and leaves the other coordinates
unchanged, apart from relabeling.
If the last positive coordinate is at $\beta<n$, the first relation
in \eqref{NB} gives a representative containing
$x_\beta\prec x_{\beta+1}$. Its contraction is normal and replaces
$(d_\beta,0)$ by $d_\beta-1$. In both cases, normalized reinsertion
recovers the original normal monomial.
\end{lemma}
\begin{Proof}
Proceed by induction in the normal grammar. At a $\prec$ root, the
first positive coordinate is in the left factor. Adding one at its
center cannot create an earlier positive coordinate unless this factor
is a leaf, in which case the first positive position is $1$.
At a $\succ$ root with nontrivial left factor, the first positive
coordinate is again there. If the left factor is a leaf, the first
positive coordinate lies in the pure right factor, unless both factors
are leaves. Induction therefore gives the required $\succ$-submonomial.
Its contraction preserves normality and removes one vertex centered
at $\alpha$. All ancestor centers there pass to the merged leaf.

For the last positive coordinate, write the normal root as $U\prec P$.
If $P$ is nontrivial, the required $\prec$-submonomial is obtained
inside this pure factor by the reflected pure argument. Suppose
$P=x_n$. If $U$ is a leaf, the assertion is immediate. If the root of
$U$ is $\prec$, apply induction inside $U$ and retain the attachment
$\prec x_n$.
Otherwise, $U=A\succ Q$ with $Q$ pure $\succ$, and $\beta=n-1$.
Move $\prec x_n$ down the right spine using only
$(A\succ B)\prec C=A\succ(B\prec C)$ until it reaches
$x_{n-1}\prec x_n$. Contracting this submonomial recovers $U$ with
a relabeled last leaf. These moves preserve centers, which proves the
asserted degree formula.
\end{Proof}

\begin{lemma}
\label{lem:two-positive}
Let $r,s\geq2$, $m\geq0$, $n=r+m+s+2$ and $h=n-2$.
Suppose the only nonsingleton positions are $\alpha=r+1$ and
$\beta=r+m+2$. Let $N_a$ be the normal monomial with
$d_\alpha=a$, $d_\beta=h+1-a$, $1\leq a\leq h$.
The only $N_a$ with a rigid proper end are those with $a=r$ or
$a=r+m+1$.
Assume that in smaller arities every sector containing the four
boundary singletons has one fibre among its normal monomials without
rigid proper ends. Then $N_h-N_1\in D_n$, except in arity six,
where the endpoints are joined by four applications of \eqref{R}.
\end{lemma}
\begin{Proof}
Write the leaves in order as
$a_1,\ldots,a_r,u,z_1,\ldots,z_m,v,b_1,\ldots,b_s$, and put
$A_i=S(a_i)$, $Z_i=S(z_i)$ and $B_i=S(b_i)$. The operated images are
\begin{align}
\phi(N_a)&=A_1\cdots A_{r-a}
S(A_{r-a+1}\cdots A_r u)Z_1\cdots Z_m vB_1\cdots B_s,
&&1\leq a\leq r,\label{two-left}\\
\phi(N_a)&=S(A_1\cdots A_r uZ_1\cdots Z_{a-r})
Z_{a-r+1}\cdots Z_m vB_1\cdots B_s,
&&r\leq a\leq r+m,\label{two-middle}\\
\phi(N_a)&=A_1\cdots A_r uZ_1\cdots Z_m
S(vB_1\cdots B_{h+1-a})B_{h+2-a}\cdots B_s,
&&r+m+1\leq a\leq h.\label{two-right}
\end{align}
Empty strings of factors are omitted.

For \eqref{two-left}, start with the right $\succ$-comb on
$a_{r-a+1},\ldots,a_r,u$. Attach the right comb on
$z_1,\ldots,z_m,v$ by $\succ$, prepend $a_1,\ldots,a_{r-a}$ as a right comb,
and then attach each $b_i$ by $\prec$.
For \eqref{two-middle}, start with the right comb on
$a_1,\ldots,a_r,u$, attach $z_1,\ldots,z_{a-r}$ by $\prec$, then
attach by $\succ$ the right comb on the remaining $z_i,v$, and finally
attach the $b_i$ by $\prec$.
For \eqref{two-right}, attach all $z_i$ by $\prec$ to the initial
right comb, then attach the left $\prec$-comb on
$v,b_1,\ldots,b_{h+1-a}$, and finally attach the remaining $b_i$
individually.
Each construction is normal and has the stated degrees. Theorem~\ref{thm:degreevector}
therefore proves the formulas.

Only $N_r$ has a rigid global initial operator, and only $N_{r+m+1}$
has a rigid global final operator. In every other case the corresponding
boundary has no nontrivial operator, or its argument ends, respectively
starts, in an $S$-factor and admits the opposite factorization. Thus
there are no other rigid ends.

Suppose $m\geq1$. Delete one middle singleton and use primes for the
sector $(r,m-1,s)$ of arity $n-1$. Let $I_p$ insert $\prec$ at its
first positive input, and let $I_s$ insert $\succ$ at its second
positive input, normalizing afterwards. The mixed relations give
\begin{align}
I_p(N'_{h-1})&=N_h,& I_p(N'_1)&=N_{r+1},\label{two-insert-p}\\
I_s(N'_1)&=N_1,& I_s(N'_{h-1})&=N_{r+m}.\label{two-insert-s}
\end{align}
For the second equality in \eqref{two-insert-p}, move the inserted
$\prec$ out through the $r-1$ outer $\succ$-vertices by
$a\succ(B\prec C)=(a\succ B)\prec C$, and then use
$(A\prec b)\succ C=A\prec(b\succ C)$. The resulting normal image is
$S(A_1\cdots A_r uZ_1)Z_2\cdots Z_m vB_1\cdots B_s$.
For the second equality in \eqref{two-insert-s}, the relevant part
before insertion is $C\prec(v\prec b_1)$. Substituting
$v\mapsto z_m\succ v$ gives
\begin{multline*}
C\prec((z_m\succ v)\prec b_1)
=C\prec(z_m\succ(v\prec b_1))=(C\prec z_m)\succ(v\prec b_1)
=((C\prec z_m)\succ v)\prec b_1.
\end{multline*}
The remaining attachments are unchanged, giving \eqref{two-middle}
with $a=r+m$. The first equality in \eqref{two-insert-p} follows by
moving the new $\prec$ through the initial right comb using the first
mixed relation. The first equality in \eqref{two-insert-s} is literal.

If $m=1$, the two images of the smaller fibre intersect at $N_{r+1}$.
If $m\geq2$, the vertex $N'_{r+1}$ has no rigid end, since the two
excluded parameters in the smaller sector are $r$ and $r+m$.
Its second-positive insertion is literal and gives
$I_s(N'_{r+1})=N_{r+1}$. Thus the images also intersect in this case.
By the hypothesis, all differences within either image belong to $D_n$.
Together the images contain $N_h,N_1$, so $N_h-N_1\in D_n$.

Now suppose $m=0$ and $s\geq3$. The monomial $N_h$ is the graft of
the right $\succ$-comb on $a_1,\ldots,a_r,u$ and a pure $\prec$
factor on $u,v,b_1,\ldots,b_s$ with vector $(s,1,0^s)$.
A reflected $F_3$ changes this proper factor to vector $(s-1,2,0^s)$.
Before reflection, the identity is
\[
c_{s-2}\succ\bigl(c_{s-1}\succ((c_s\succ x)\succ y)\bigr)
=c_{s-2}\succ\bigl((c_{s-1}\succ(c_s\succ x))\succ y\bigr).
\]
Prepending $c_1,\ldots,c_{s-3}$ as a right comb gives vectors
$(0^s,1,s)$ and $(0^s,2,s-1)$; for $s=3$, this prefix is absent.
After reflection, the whole graft becomes $N_{h-1}$ modulo $D_n$.
Both $N_{h-1}$ and $N_1$ have the last $\prec$-submonomial from
Lemma~\ref{lem:extreme-cherry}. Contracting it gives the two
concentrated monomials in the smaller sector $(r,0,s-1)$.
The hypothesis and reinsertion give $N_{h-1}-N_1\in D_n$.
Reflection treats $s=2$, $r\geq3$.

Finally, $r=s=2$, $m=0$ has arity six. Retain the four singleton
intervals and change the nonsingleton interval successively through
\[
[2,3],\quad[1,2],\quad[2,5],\quad[5,6],\quad[4,5].
\]
These are four valid moves, with singleton pivots $2,2,5,5$.
Their endpoints are $N_1,N_4$, and every substituted block is nonempty.
\end{Proof}

\begin{theorem}
\label{thm:strong}
For every $\{1,2,n-1,n\}\subseteq M\subsetneq\{1,\ldots,n\}$,
the normal monomials without rigid proper ends form one operated fibre.
For $n\geq8$, all $A(M,k)$ are equal modulo $D_n$.
Consequently, the kernel on this part of the sector is spanned modulo
$D_n$ by the arity-$n$ corner identities $J,H,\rho J$; for $|M|>4$,
it is contained in $D_n$.
\end{theorem}
\begin{Proof}
Proceed by induction on $n$. There is no such sector below arity five.
By Proposition~\ref{prop:concentrated}, every monomial under
consideration reduces to some $A(M,k)$. If there is one nonsingleton
position, this form is unique; this includes arity five.
If there are two, Lemma~\ref{lem:two-positive} connects the two forms
using the induction hypothesis and supplies the arity-six case.

Suppose there are at least three nonsingleton positions, with first
$\alpha$ and last $\beta$. For $k\ne\alpha$, the first positive
degree of $A(M,k)$ is one. Contract the corresponding $\succ$
submonomial from Lemma~\ref{lem:extreme-cherry}. Both boundary run
lengths and the singleton count are unchanged, and the new vector is
concentrated at the image of $k$. Its mobile graft has no rigid ends.
Induction connects these smaller forms. Reinsertion shows that all
$A(M,k)$ with $k\ne\alpha$ are equal modulo $D_n$.
Similarly, contraction at the last positive position and normalized
reinsertion identify all forms with $k\ne\beta$.
These two sets intersect because a third nonsingleton position exists.
Their union contains every $A(M,k)$.

This proves the fibre assertion and the required congruence of the
concentrated forms. The last assertion now follows from
Proposition~\ref{prop:concentrated}.
\end{Proof}

\noindent\textbf{Completion of the upper bound.}
\label{sec:upper-final}

Let $\mathcal V_n$ be the span of the arity-$n$ pure consequences of
$F_r,G_r,K_6,E_6$ and their reflections. Let $\mathcal U_n$ be the
sum of $\mathcal V_n$ and the spans of all arity-$n$ instances of
\begin{align*}
B(U,z,V)&=U\prec(z\prec V)-(U\succ z)\prec V,
&&\LL(\phi(U)),\ \RR(\phi(V)),\\
E(a,b,c,d;V)&=(a\succ(b\succ c))\prec(d\prec V)
-(a\succ((b\succ c)\succ d))\prec V,
&&\RR(\phi(V)),
\end{align*}
their reflections, and $F_m,\rho F_m$ with arbitrary nonempty binary
arguments. By \eqref{generalB}, \eqref{generalE} and
Theorem~\ref{thm:validity}, $\mathcal U_n\subseteq K_n$.
The stated families $F,G,J,H$ belong to this space.

\begin{theorem}
\label{thm:upper}
Over any field and for every $n\geq8$,
$$K_n\subseteq D_n+\Span\{F_{n-2},\rho F_{n-2},G_{n-3},\rho G_{n-3},
J_{n-4},\rho J_{n-4},H_{p,n-p-1}:4\leq p\leq n-5\}.$$
In particular, $K_n=D_n+\mathcal U_n$ and $\dim Q_n\leq n-2$.
\end{theorem}
\begin{Proof}
By Lemma~\ref{lem:fibres}, consider a binomial $T-T'$ in one fibre,
with singleton set $M$. If either endpoint has a nontrivial rigid
proper end, Lemma~\ref{lem:old-rigid} gives $T-T'\in D_n$.
Assume there are no such ends.

If $1,n\notin M$, apply Theorem~\ref{thm:interior}.
If exactly one of $1,n$ lies in $M$, apply
Proposition~\ref{prop:one-boundary}. A boundary run of length one gives
a difference in $D_n$; otherwise the class is pure. By pure
completeness, only $F_{n-2},G_{n-3}$ or their reflections remain in
arity $n$.
If $1,n\in M$ and either run has length one, apply
Lemma~\ref{lem:short-run}. Finally, if
$\{1,2,n-1,n\}\subseteq M$, Theorem~\ref{thm:strong} reduces the
class to the corner span.

This proves the inclusion. The displayed generators belong to
$\mathcal U_n\subseteq K_n$, so it also gives $K_n=D_n+\mathcal U_n$.
There are six classes from $F,G,J$ and their reflections and $n-8$
classes from $H$, with the latter absent for $n=8$. Thus
$\dim Q_n\leq n-2$.
\end{Proof}

\section{Independence and the main theorem}

We prove the lower bound independently of Theorem~\ref{thm:upper}.
Let $\iota_i^\prec$ and $\iota_i^\succ$ insert
$x_i\prec x_{i+1}$ and $x_i\succ x_{i+1}$ at the $i$th leaf,
with the remaining labels shifted. For $\mu\in\{\prec,\succ\}$,
put $\lambda_\mu(f)=x_1\mu f(x_2,\ldots,x_n)$ and
$r_\mu(f)=f(x_1,\ldots,x_{n-1})\mu x_n$.

\begin{lemma}
\label{lem:one-step-old}
For every $n\geq2$,
\[
D_n=\sum_{\mu\in\{\prec,\succ\}}
\left(\sum_{i=1}^{n-1}\iota_i^\mu(K_{n-1})
+\lambda_\mu(K_{n-1})+r_\mu(K_{n-1})\right).
\]
\end{lemma}
\begin{Proof}
The right side is contained in $D_n$. Conversely, remove the last
binary vertex added in constructing an arity-$n$ consequence.
The remaining relation has arity $n-1$ and belongs to $K_{n-1}$.
Restoring the vertex gives one of the stated insertions or extensions.
\end{Proof}

\begin{lemma}
\label{lem:isolated-G-source}
For $m\geq5$, the normal monomial in the sector $\{1,2,m-1\}$ with
excess positions $3,m$ is a singleton $S$-fibre.
\end{lemma}
\begin{Proof}
Its interval family is
$\{[1,1],[2,2],[m-1,m-1]\}\cup\{[1,j]:3\leq j\leq m-2\}$.
An interval starting at $1$ cannot be a pivot. The pivot $[2,2]$
requires $[1,2]$ or a larger interval starting at $2$, neither of
which is present. The pivot $[m-1,m-1]$ requires a larger interval
starting or ending there, which is also absent. Thus \eqref{R}
cannot be applied in either direction.
\end{Proof}

\begin{proposition}
\label{prop:FG-independent}
For every $n\geq8$, the classes of
$F_{n-2},\rho F_{n-2},G_{n-3},\rho G_{n-3}$ are nonzero and linearly
independent in $Q_n$.
\end{proposition}
\begin{Proof}
Their singleton sets are, respectively,
$\{1,2,3\}$, $\{n-2,n-1,n\}$, $\{1,2,n-1\}$ and $\{2,n-1,n\}$.
By \eqref{sector-decomp} and reflection, it suffices to prove
nonvanishing for $F_{n-2}$ and $G_{n-3}$.

The second endpoint of $F_{n-2}$ is
\[
f_n=x_1\succ L_{n-3}\bigl(x_2\succ(x_3\succ x_4),x_5,\ldots,x_n\bigr).
\]
The relations of $\NB$ preserve the number of $\prec$-vertices.
Hence a one-step consequence can have nonzero coefficient at the pure
monomial $f_n$ only if its source is pure and the added operation is
$\succ$. The only two-leaf submonomial of $f_n$ is $x_3\succ x_4$;
contracting it gives the sector $\{1,2\}$, whose kernel is zero.
The only outer contraction gives the same sector. By
Lemma~\ref{lem:one-step-old}, the coefficient at $f_n$ therefore
vanishes on $D_n$, whereas it is $\pm1$ on $F_{n-2}$.

For $G_{n-3}$, use the second endpoint
$g_n=T_{n-2}(x_1,\ldots,x_{n-2})\succ(x_{n-1}\succ x_n)$.
Its two-leaf submonomials are $x_2\succ x_3$ and
$x_{n-1}\succ x_n$, and neither root factor is a leaf.
Contracting the latter gives the zero-kernel pair sector. Contracting
the former gives the singleton fibre in
Lemma~\ref{lem:isolated-G-source}. Thus the coefficient at $g_n$
also vanishes on $D_n$ and is $\pm1$ on $G_{n-3}$.
\end{Proof}

Retain $K=n-4$ and the corner notation from Section~\ref{sec:corner}.

\begin{lemma}
\label{lem:corner-sources}
Put $m=n-1$. A one-leaf insertion with target singleton set
$\{1,2,n-1,n\}$ has source in one of the sectors
\[
\{1,2,m-1,m\},\qquad\{1,2,m\},\qquad\{1,m-1,m\}.
\]
The first sector enters through the two maps \eqref{corner-embeddings}.
An outer extension comes from one of the two triple sectors.
\end{lemma}
\begin{Proof}
A $\prec$-insertion at input $i$ makes $i$ a nonsingleton and
$i+1$ a singleton. The required target set therefore forces
$i=n-2$. After contraction, the source singleton set is
$\{1,2,m\}$, with $m-1$ possibly added. A $\succ$-insertion makes
$i$ a singleton and $i+1$ a nonsingleton, so $i=2$. Its source set
is $\{1,m-1,m\}$, with $2$ possibly added. For corner sources,
\eqref{corner-excess} gives precisely the maps
\eqref{corner-embeddings}.

Of the four outer extensions, $x_1\prec T$ has a nonsingleton first
leaf and $T\succ x_n$ has a nonsingleton last leaf. The other two
extensions add one boundary singleton and hence have source in
the indicated triple sectors.
\end{Proof}

\begin{proposition}
\label{prop:corner-independent}
For every $n\geq8$, the $n-6$ classes
\[
J_{n-4},\quad H_{p,n-p-1}\ (4\leq p\leq n-5),\quad\rho J_{n-4}
\]
are linearly independent modulo $D_n$.
\end{proposition}
\begin{Proof}
Identify the corner basis with $\Delta_K$, and let
$\pi_K\colon\Span\Delta_K\to\Span\mathcal E_K$ be coordinate
projection. By Lemma~\ref{lem:one-step-old}, it suffices to consider
one-step consequences. The two triple source kernels vanish by
Proposition~\ref{cor:triple-zero-kernel}.
Under the unshifted embedding, every preimage of a point of
$\mathcal E_K$ has terminal coordinate $K$ and is locked. Under the
shifted embedding, every such preimage has $a=b\geq2$ and is again
locked. Lemma~\ref{lem:corner-locked} therefore gives
$\pi_K(D_n\cap K_n(M_0))=0$.

For $n=8$, the images of $J_4$ and $\rho J_4$ are distinct coordinates.
Suppose $n\geq9$. Write $e_0=(1,1,K)$,
$e_t=(t,t+1,K)$ for $2\leq t\leq K-3$, and
$e_\infty=(K-1,K,K+1)$, also using these symbols for their basis
vectors. By \eqref{Jcorner}--\eqref{rhoJcorner}, up to signs,
\[
\pi_K(J_{n-4})=e_0,\qquad
\pi_K(H_{4,n-5})=e_2,\qquad
\pi_K(\rho J_{n-4})=e_\infty,
\]
and $\pi_K(H_{p,n-p-1})=e_{p-3}-e_{p-2}$ for
$5\leq p\leq K-1$. In the basis
$e_0,e_2,\ldots,e_{K-3},e_\infty$, the two extreme blocks are
$[\pm1]$ and the middle block is bidiagonal with diagonal entries
$\pm1$. Its determinant is $\pm1$, proving independence over any field.
\end{Proof}

The reflected $H$-identities give no additional classes. For
$p\geq5$, put $W=T_{q+1}$ and write $Q_p=u\prec Q_{p-1}$ on the
remaining consecutive variables. Then
\[
\rho H_{p,q}=W\succ Q_p-W\prec Q_p
=(W\succ u)\prec Q_{p-1}-W\prec(u\prec Q_{p-1})
=-H_{q+1,p-1}
\]
already modulo \eqref{NB}.

For $p=4$, put $n=q+5$, $K=n-4$ and $Z_j=T_j\prec Q_{n-j}$.
Normalizing the reflected endpoints gives
$\rho H_{4,q}=A-Z_{q+1}$, where
$A=(T_{q+2}\prec(x_{q+3}\prec x_{q+4}))\prec x_{q+5}$.
The coordinates of $A$ and $Z_4$ are $(K-2,K-1,K)$ and $(1,2,K)$,
respectively. Both lie in $\mathcal O_K$, so $A-Z_4\in D_n$.
Since $H_{j,n-j-1}=Z_j-Z_{j+1}$, telescoping gives
\[
\rho H_{4,q}=\sum_{j=4}^{q}H_{j,q+4-j}\qquad\text{in }Q_{q+5}.
\]

\begin{theorem}
\label{thm:classification}
For every $n\geq8$, $\dim Q_n=n-2$, and a basis is given by the classes
of
$$F_{n-2},\ \rho F_{n-2},\ G_{n-3},\ \rho G_{n-3},\
J_{n-4},\ \rho J_{n-4},H_{p,n-p-1},\qquad4\leq p\leq n-5.$$
\end{theorem}
\begin{Proof}
Theorem~\ref{thm:upper} gives $\dim Q_n\leq n-2$.
Proposition~\ref{prop:FG-independent} gives four independent classes
in distinct three-singleton sectors, and
Proposition~\ref{prop:corner-independent} gives $n-6$ independent
classes in the four-singleton corner. By \eqref{sector-decomp},
these classes are independent together. Thus
$\dim Q_n\geq4+(n-6)=n-2$.
\end{Proof}

\noindent\textbf{The initial generators.}

Denote the identities \eqref{chi} and \eqref{Cfive} by $\chi_4$ and
$\chi_5$, respectively. We also use
\begin{align*}
\chi_6:\quad
\bigl((a\succ(b\succ c))\prec(d\prec e)\bigr)\prec f
&=\bigl(T_4(a,b,c,d)\prec e\bigr)\prec f,\\
\chi_7:\quad
(a\succ(b\succ c))\prec\bigl((d\prec(e\prec f))\prec g\bigr)
&=\bigl(T_4(a,b,c,d)\prec(e\prec f)\bigr)\prec g.
\end{align*}
As before, we identify an identity with the difference of its two sides.

\begin{proposition}
\label{prop:initial-generators}
We have $K_1=K_2=K_3=0$. For $4\leq n\leq7$, a basis of $Q_n$
is given by the classes of the following sets:
\begin{align*}
\mathcal G_4&=\{\chi_4\},\\
\mathcal G_5&=\{F_3,\rho F_3,\chi_5,\rho\chi_5\},\\
\mathcal G_6&=\{F_4,\rho F_4,G_3,\rho G_3,K_6,\rho K_6,\chi_6\},\\
\mathcal G_7&=\{F_5,\rho F_5,G_4,\rho G_4,\chi_7,\rho\chi_7\}.
\end{align*}
In particular, $\dim Q_4=1$, $\dim Q_5=4$, $\dim Q_6=7$ and
$\dim Q_7=6$, over every field.
\end{proposition}

\begin{Proof}
The validity of $\chi_6$ follows from the arity-six case of
Lemma~\ref{lem:two-positive}. The endpoints of $\chi_7$ have corner
coordinates $(1,1,3)$ and $(1,2,3)$, so its validity follows from
the paths in Section~\ref{sec:corner}. The other identities were
proved above.

Put $I_3=0$ and let $I_n$ be the operadic ideal of $\NB$ generated
by $\mathcal G_4,\ldots,\mathcal G_n$. For $4\leq n\leq7$, write
\[
c_n=\dim(\NB/I_{n-1})(n),\qquad
b_n=\dim(\NB/I_n)(n),\qquad
f_n=\dim\operatorname{im}\phi_n.
\]
We compute $c_n$ and $b_n$ using the program of Dotsenko and
Heijltjes~\cite{DotsHij}, in the nonsymmetric case. The input
consists of \eqref{NB} and the generators of $I_{n-1}$ or $I_n$,
respectively. We compute Gr\"obner bases through arity $n$ and
count the normal monomials of that arity. All critical pairs of
arity at most $n$ reduce to zero. Moreover, every resulting rule
equates two monomials, so these reductions and the resulting
dimensions are valid over every field.

To compute $f_n$, close $\phi(\Ncal_n)$ under both directions
of \eqref{R}, including applications inside $S$-arguments.
This gives a finite graph, since every word has $n$ ordered
variables and $n-1$ occurrences of $S$. By
Lemma~\ref{lem:fibres}, $f_n$ is the number of its connected
components. In arities at most three, distinct normal monomials
belong to distinct components, giving $K_1=K_2=K_3=0$.
For the remaining arities, the computations give
\begin{center}
\begin{tabular}{c|rrrrrr}
$n$ & $|\Ncal_n|$ & $c_n$ & $b_n$ & $f_n$
& $\dim K_n$ & $c_n-b_n$\\
\hline
4 & 20  & 20  & 19  & 19  & 1   & 1\\
5 & 70  & 60  & 56  & 56  & 14  & 4\\
6 & 252 & 158 & 151 & 151 & 101 & 7\\
7 & 924 & 391 & 385 & 385 & 539 & 6
\end{tabular}
\end{center}

Since $I_n\subseteq K$ and $b_n=f_n$, we have $I_n(n)=K_n$
for $4\leq n\leq7$. Proceeding by induction on $n$, we obtain
$I_{n-1}(n)=D_n$. Hence the classes of $\mathcal G_n$ span $Q_n$,
and $\dim Q_n=c_n-b_n=|\mathcal G_n|$. Therefore these classes
form a basis.
\end{Proof}

Put $\mathcal G_{\leq7}=\mathcal G_4\cup\mathcal G_5\cup
\mathcal G_6\cup\mathcal G_7$.
\begin{theorem}
\label{thm:generators}
The kernel $K=\ker\phi$ has a minimal homogeneous generating set
consisting of the eighteen identities in
Proposition~\ref{prop:initial-generators} and
\begin{gather*}
F_r,\rho F_r\ (r\geq6),\qquad G_r,\rho G_r\ (r\geq5),\\
J_r,\rho J_r\ (r\geq4),\qquad H_{p,q}\ (p,q\geq4).
\end{gather*}
\end{theorem}
\begin{Proof}
Proceed by induction on arity. Proposition~\ref{prop:initial-generators} gives generation through
arity seven. Once the smaller components are generated, their arity-$n$
consequences span $D_n$, and Theorem~\ref{thm:classification} gives
a basis of $K_n/D_n$. This proves generation.
If a generator is omitted, its class cannot be expressed through the
remaining basis classes or through consequences of smaller arity.
Larger arities cannot generate a smaller component, since there are
no nullary operations. Hence the generating set is minimal.
\end{Proof}

Since $Q_n\ne0$ for every $n\geq8$, the operadic ideal $K$ is not
finitely generated.

\begin{remark}
If the associative multiplication is commutative, then
$a\succ b=b\prec a$, and $a\circ b=aS(b)$ satisfies
\[
(a\circ b)\circ c=(a\circ c)\circ b,\qquad
a\circ(b\circ c)=b\circ(a\circ c).
\]
Every identity of the derived operation in the commutative case
follows from these two identities.

Indeed, let $A=\kk[y_i,z_i\mid i\geq1]_+$ be the algebra of
commutative polynomials without constant term. Define $S$ on monomials
by $S(y_i)=z_i$, $S(z_i)=z_i$, and $S(m)=m$ for $\deg m\geq2$.
Since $S|_{A^2}=\id$, we have
$aS(S(b)c)=aS(b)c=S(aS(b))c$, so \eqref{R} holds.

Let $B$ be the subalgebra generated by the $y_i$ under $u\circ v=uS(v)$.
For mixed monomials $u=Y^\alpha Z^\beta$ and $v=Y^\gamma Z^\delta$,
where $|\alpha|,|\beta|,|\gamma|,|\delta|>0$, the multiplication is
\[
y_i\circ y_j=y_i z_j,\qquad y_i\circ u=y_i u,\qquad
u\circ y_j=uz_j,\qquad u\circ v=uv.
\]
Left multiplication by the $y_i$ adds $y$-factors, and right
multiplication adds $z$-factors. Thus $B$ has basis $y_i$ and
$Y^\alpha Z^\beta$ with $|\alpha|,|\beta|>0$. These are precisely
the basis and multiplication of the free bicommutative algebra over
an arbitrary field \cite[Theorem~2.1]{DIT}; see also
\cite[p.~173]{DrenskyBicom}. Hence $x_i\mapsto y_i$ embeds it into
$(A,\circ)$.

Thus the suboperad generated by $x_1S(x_2)$ in the commutative case
is isomorphic to the bicommutative operad. There are no additional
special identities, even under the assumptions $S^2=S$ and
$S|_{A^2}=\id$.
\end{remark}


\begin{thebibliography}{99}

\bibitem{Aguiar}
M. Aguiar,
\href{https://doi.org/10.1023/A:1010818119040}{Pre-Poisson algebras},
\textit{Lett. Math. Phys.} \textbf{54} (2000), 263--277.

\bibitem{Bergman}
G. M. Bergman,
\href{https://doi.org/10.1016/0001-8708(78)90010-5}{The diamond lemma for ring theory},
\textit{Adv. Math.} \textbf{29} (1978), 178--218.

\bibitem{BCQ}
L. A. Bokut, Y. Chen, J. Qiu,
\href{https://arxiv.org/abs/0805.0640}{Gr\"obner--Shirshov bases for associative algebras with multiple operators and free Rota--Baxter algebras},
\textit{J. Pure Appl. Algebra} \textbf{214} (2010), 89--100.

\bibitem{BD}
M. R. Bremner, V. Dotsenko,
\href{https://doi.org/10.1201/b20061}{\textit{Algebraic Operads: An Algorithmic Companion}},
CRC Press, Boca Raton, 2016.

\bibitem{DSBasis}
A. Dauletiyarova, B. K. Sartayev,
\href{https://doi.org/10.1142/S0219498825502925}{Basis of the free noncommutative Novikov algebra},
\textit{J. Algebra Appl.} \textbf{24} (2025), no. 12, 2550292.

\bibitem{DotsHij}
V.~Dotsenko, W.~Heijltjes,
Gr\"obner bases for operads,
Haskell software.
\url{https://vladimirdotsenko.github.io/}.

\bibitem{DK}
V. Dotsenko, A. Khoroshkin,
\href{https://doi.org/10.1215/00127094-2010-026}{Gr\"obner bases for operads},
\textit{Duke Math. J.} \textbf{153} (2010), 363--396.

\bibitem{DrenskyBicom}
V. Drensky,
\href{https://www.math.bas.bg/serdica/2019/2019-167-188.pdf}{Varieties of bicommutative algebras},
\textit{Serdica Math. J.} \textbf{45} (2019), 167--188.

\bibitem{DZ}
V. Drensky, B. K. Zhakhayev,
\href{https://doi.org/10.1016/j.jalgebra.2017.12.012}{Noetherianity and Specht problem for varieties of bicommutative algebras},
\textit{J. Algebra} \textbf{499} (2018), 570--582.

\bibitem{DI}
A. S. Dzhumadil'daev, N. A. Ismailov,
\href{https://arxiv.org/abs/1711.04300}{Polynomial identities of bicommutative algebras, Lie and Jordan elements},
\textit{Comm. Algebra} \textbf{46} (2018), 5241--5251.

\bibitem{DIT}
A. S. Dzhumadil'daev, N. A. Ismailov, K. M. Tulenbaev,
\href{https://www.math.bas.bg/serdica/2011/2011-025-044.pdf}{Free bicommutative algebras},
\textit{Serdica Math. J.} \textbf{37} (2011), 25--44.

\bibitem{DL}
A. S. Dzhumadil'daev, C. L\"ofwall,
\href{https://www.kurims.kyoto-u.ac.jp/EMIS/journals/HHA/volumes/2002/n2a7/v4n2a7.pdf}{Trees, free right-symmetric algebras, free Novikov algebras and identities},
\textit{Homology Homotopy Appl.} \textbf{4} (2002), no. 2, part 1, 165--190.

\bibitem{EFGuo}
K. Ebrahimi-Fard, L. Guo,
\href{https://arxiv.org/abs/math/0503647}{Rota--Baxter algebras and dendriform algebras},
\textit{J. Pure Appl. Algebra} \textbf{212} (2008), 320--339.

\bibitem{GK}
V. Ginzburg, M. Kapranov,
\href{https://doi.org/10.1215/S0012-7094-94-07608-4}{Koszul duality for operads},
\textit{Duke Math. J.} \textbf{76} (1994), 203--272.

\bibitem{Guo}
L. Guo,
\textit{An Introduction to Rota--Baxter Algebra},
Surveys of Modern Mathematics, vol. 4,
International Press, Somerville, MA, and Higher Education Press, Beijing, 2012.

\bibitem{KMS}
P. Kolesnikov, F. Mashurov, B. Sartayev,
\href{https://doi.org/10.3842/SIGMA.2024.017}{On pre-Novikov algebras and derived Zinbiel variety},
\textit{SIGMA} \textbf{20} (2024), 017, 15 pp.

\bibitem{KSPerm}
P. S. Kolesnikov, B. K. Sartayev,
\href{https://doi.org/10.1080/00927872.2022.2028798}{On the embedding of left-symmetric algebras into differential Perm-algebras},
\textit{Comm. Algebra} \textbf{50} (2022), 3246--3260.

\bibitem{KSGD}
P. S. Kolesnikov, B. K. Sartayev,
\href{https://doi.org/10.1080/10586458.2022.2041134}{On the special identities of Gelfand--Dorfman algebras},
\textit{Exp. Math.} \textbf{33} (2024), 165--174.

\bibitem{KSO}
P. S. Kolesnikov, B. Sartayev, A. Orazgaliev,
\href{https://doi.org/10.1016/j.jalgebra.2019.07.034}{Gelfand--Dorfman algebras, derived identities, and the Manin product of operads},
\textit{J. Algebra} \textbf{539} (2019), 260--284.

\bibitem{LM}
Z. Liu, P. Malbos,
\href{https://arxiv.org/abs/2502.16304}{Polygraphic resolutions for operated algebras},
arXiv:2502.16304, 2025.

\bibitem{Loday}
J.-L. Loday,
\href{https://arxiv.org/abs/math/0102053}{Dialgebras},
in \textit{Dialgebras and Related Operads}, Lecture Notes in Mathematics,
vol. 1763, Springer, Berlin, 2001, 7--66.

\bibitem{LV}
J.-L. Loday, B. Vallette,
\href{https://doi.org/10.1007/978-3-642-30362-3}{\textit{Algebraic Operads}},
Grundlehren der mathematischen Wissenschaften, vol. 346,
Springer, Heidelberg, 2012.

\bibitem{MS}
F. A. Mashurov, B. K. Sartayev,
\href{https://arxiv.org/abs/2605.01007}{Nonsymmetric versions of binary quadratic operads},
arXiv:2605.01007, 2026.

\bibitem{SKNov}
B. Sartayev, P. Kolesnikov,
\href{https://doi.org/10.1007/s40879-023-00632-1}{Noncommutative Novikov algebras},
\textit{Eur. J. Math.} \textbf{9} (2023), 35.

\bibitem{Vallette}
B. Vallette,
\href{https://arxiv.org/abs/math/0609002}{Manin products, Koszul duality, Loday algebras and Deligne conjecture},
\textit{J. Reine Angew. Math.} \textbf{620} (2008), 105--164.

\end{thebibliography}
\end{document}